%% file: DESW-SlopesSegments.tex
\documentclass[a4paper,11pt]{article}

\input{font}
\usepackage{amsfonts,amsmath,graphicx,calc,enumerate,url}
\usepackage[sort&compress,numbers]{natbib}
\usepackage[lmargin=32mm,rmargin=32mm,tmargin=32mm,bmargin=32mm]{geometry}

\input{myTheorems}

\input{myLabels}

\input{myCommands}

\newcommand{\pred}[1]{\ensuremath{P{(#1)}}}

\begin{document}

\date{Submitted: July 28, 2005; Revised: \today}


\title{\textbf{Drawings of Planar Graphs with\\ Few Slopes and Segments}\,\thanks{A preliminary version of this paper was published as: ``Really straight graph drawings.'' \emph{Proceedings of the 12th International Symposium on Graph Drawing} (GD '04), \emph{Lecture Notes in Computer Science} 3383:122--132, Springer, 2004.}\\[2ex] {\emph{\large Dedicated to Godfried Toussaint on his 60th birthday.}}\\[2ex]}

\FancyFootnotes

\author{Vida Dujmovi{\'c} \footnotemark[2] \and 
David Eppstein \footnotemark[3] \and 
Matthew Suderman \footnotemark[4] \and 
David R. Wood \footnotemark[5]}

\footnotetext[2]{Department of Mathematics and Statistics, McGill University,
Montr{\'e}al, Qu\'ebec, Canada (\texttt{vida@cs.mcgill.ca}). Supported by NSERC.}

\footnotetext[3]{Department of Computer Science, University of California, Irvine, California, U.S.A. (\texttt{eppstein@ics.uci.edu}).}

\footnotetext[4]{McGill Centre for Bioinformatics, School of Computer Science, McGill University, Montr{\'e}al, Qu\'ebec, Canada (\texttt{msuder@cs.mcgill.ca}). Supported by NSERC.}

\footnotetext[5]{Departament de Matem{\`a}tica Aplicada II, Universitat Polit{\`e}cnica de Catalunya, Barcelona, Catalunya, Spain (\texttt{david.wood@upc.edu}). Supported by a Marie Curie 
Fellowship of the European Community under contract 023865, and by the 
projects MCYT-FEDER BFM2003-00368 and Gen.\ Cat 2001SGR00224.}

\maketitle

\begin{abstract}  
We study straight-line drawings of planar graphs with few segments and few slopes. Optimal results are obtained for all trees. Tight bounds are obtained for outerplanar graphs, $2$-trees, and planar $3$-trees.  We prove that every $3$-connected plane graph on $n$ vertices has a plane drawing with at most $\frac{5}{2}n$ segments and at most $2n$ slopes.  We prove that every cubic $3$-connected plane graph has a plane drawing with three slopes (and three bends on the outerface). In a companion paper, drawings of non-planar graphs with few slopes are also considered. 
\end{abstract}

\FancyFootnotesOff

\newpage

\section{Introduction}
\seclabel{Introduction}

A common requirement for an aesthetically pleasing drawing of a graph is that the edges are straight. This paper studies the following additional requirements of straight-line graph drawings:

\begin{enumerate}
\item minimise the number of segments in the drawing, and 
\item minimise the number of distinct edge slopes in the drawing.
\end{enumerate}

First we formalise these notions. Consider a mapping of the vertices of a graph to distinct points in the plane. Now represent each edge by the closed line segment between its endpoints. Such a mapping is a (\emph{straight-line}) \emph{drawing} if the only vertices that each edge intersects are its own endpoints.  By a \emph{segment} in a drawing, we mean a maximal set of edges that form a line segment.  The \emph{slope} of a line $L$ is the angle swept from the X-axis in an anticlockwise direction to $L$ (and is thus in $[0,\pi)$). The \emph{slope} of an edge or segment is the slope of the line that extends it. Of course two edges have the same slope if and only if they are parallel. A \emph{crossing} in a drawing is a pair of edges that intersect at some point other than a common endpoint. A drawing is \emph{plane} if it has no crossings. A \emph{plane graph} is a planar graph with a fixed combinatorial embedding and a specified outerface. We emphasise that a plane drawing of a plane graph must preserve the embedding and outerface. That every plane graph has a plane  drawing is a famous result independently due to  \citet{Wagner36} and \citet{Fary48}. 

In this paper we prove lower and upper bounds on the minimum number of segments and slopes in (plane) drawings of graphs. In a companion paper \citep{DSW-Slopes}, we consider drawings of non-planar graphs with few slopes. A summary of our results is given in \tabref{Summary}. A number of comments are in order when considering these results: 
\begin{itemize}
\item The minimum number of slopes in a drawing of (plane) graph $G$ is at most the minimum number of segments in a drawing of $G$.  
\item Upper bounds for plane graphs are stronger than for planar graphs, since for planar graphs one has the freedom to choose the embedding and outerface. On the other hand, lower bounds for planar graphs are stronger than for plane graphs. 
\item Deleting an edge in a drawing cannot increase the number of slopes, whereas it can increase the number of segments. Thus, the upper bounds for slopes are applicable to all subgraphs of the mentioned graph families, unlike the upper bounds for segments. 
\end{itemize}

\begin{table}[!ht]
\caption{\tablabel{Summary}Summary of results (ignoring additive constants).
Here $n$ is the number of vertices, $\eta$ is the number of vertices of odd degree, and $\Delta$ is the maximum degree. The lower bounds are existential, except for trees, for which the lower bounds are universal.}
\begin{center}
\begin{tabular}{lcccc}
\hline
graph family			& \multicolumn{2}{c}{\# segments}	& \multicolumn{2}{c}{\# slopes}\\
						& $\geq$	& $\leq$	& $\geq$	& $\leq$\\\hline
trees					& $\frac{\eta}{2}$	& $\frac{\eta}{2}$	& $\ceil{\frac{\Delta}{2}}$	& $\ceil{\frac{\Delta}{2}}$\\	
maximal outerplanar		& $n$		& $n$		& -			& $n$ 	\\
plane $2$-trees			& $2n$		& $2n$		& $2n$		& $2n$	\\
plane $3$-trees			& $2n$		& $2n$		& $2n$		& $2n$	\\
plane $2$-connected		& $\frac{5}{2}n$	& -			& $2n$		& -		\\	
planar $2$-connected	& $2n$		& -			& $n$		& -		\\
plane $3$-connected		& $2n$		& $\frac{5}{2}n$	& $2n$		& $2n$	\\	
planar $3$-connected	& $2n$		& $\frac{5}{2}n$	& $n$		& $2n$	\\
plane $3$-connected cubic	& 	-	& $n+2$		& $3$		& $3$	\\
\hline
\end{tabular}
\end{center}
\end{table}

The paper is organised as follows. In \secref{TwoThree} we consider drawings with two or three slopes, and conclude that it is \NP-complete to determine whether a graph has a plane drawing on two slopes. 

\secref{PlanarTreewidth} studies plane drawings of graphs with small treewidth. In particular, we consider trees, outerplanar graphs, $2$-trees, and planar $3$-trees. For any tree, we construct a plane drawing with the minimum number of segments and the minimum number of slopes.  For outerplanar graphs, $2$-trees, and planar $3$-trees, we determine bounds on the minimum number of segments and slopes that are tight in the worst-case.

\secref{ThreeConnected} studies plane drawings of $3$-connected plane and
planar graphs. In the case of slope-minimisation for plane graphs we obtain a
bound that is tight in the worst case. However, our lower bound examples have linear maximum degree. We drastically improve the upper bound in the case of cubic graphs. We prove that every $3$-connected plane cubic graph has a plane drawing with three slopes,  except for three edges on the outerface that have their own slope. As a corollary we prove that every $3$-connected plane cubic graph has a plane `drawing' with three slopes and three bends on the
outerface.

We now review some related work from the literature.

Plane orthogonal drawings with two slopes (and few bends) have been extensively studied \citep{RNN-CG02, RNN-JAlg00, RNN-JGAA99, RNN-CG98, Thomassen84, BS-Algo88, BS-Networks87, Ungar-JLMS53, RNG-JAlg04}.  For example, \citet{Ungar-JLMS53} proved  that every cyclically $4$-edge-connected plane cubic graph has a plane drawing with two slopes and four bends on the
outerface. Thus our result for $3$-connected plane cubic graphs (\corref{CubicBends})  nicely complements this theorem of Ungar.




Contact and intersection graphs of segments in the plane with few slopes is an interesting line of research. The \emph{intersection graph} of a set of segments has one vertex for each segment, and two vertices are adjacent if and only if the corresponding segments have a non-empty intersection. \citet{HNZ-DM91} proved that every bipartite planar graph is the intersection graph of some set of horizontal and vertical segments. A \emph{contact graph} is an intersection graph of segments for which no two segments have an interior point in common. Strengthening the above result, \citet{dFdMP95} (and later, \citet{CKU-IPL98}) proved that every bipartite planar graph is a contact graph of some set of horizontal and vertical segments. Similarly, \citet{CCDMN-JGAA02} proved that every triangle-free planar graph is a contact graph of some set of segments with only three distinct slopes. It is an open problem whether every planar graph is the  intersection graph of a set of segments in the plane; see \citep{dFdM-GD04,OdMdF99} for the most recent results. It is even possible that every $k$-colourable planar graph ($k\leq4$) is the intersection graph of some set of segments using only $k$ distinct slopes.


\mySubSection{Definitions}{Definitions}

We consider  undirected, finite, and simple graphs $G$ with  vertex set $V(G)$ and edge set $E(G)$. The  number of vertices and edges of $G$ are respectively denoted by $n=|V(G)|$ and $m=|E(G)|$. The maximum degree of $G$ is denoted by $\Delta(G)$. 

For all $S\subseteq V(G)$, the (\emph{vertex-}) \emph{induced} subgraph $G[S]$ has vertex set $S$ and edge set $\{vw\in E(G):v,w\in S\}$.  For all $S\subseteq V(G)$, let $G\setminus S$ be the subgraph $G[V(G)\setminus S]$. For all $v\in V(G)$, let $G\setminus v=G\setminus\{v\}$. For all $A,B\subseteq V(G)$, let $G[A,B]$ be the bipartite subgraph of $G$ with vertex set $A\cup B$ and edge set $\{vw\in E(G):v\in A\setminus B,w\in B\setminus A\}$. 

For all $S\subseteq E(G)$, the (\emph{edge-}) \emph{induced} subgraph $G[S]$ has vertex set $\{v\in V(G):\exists\,vw\in S\}$ and edge set $S$.  For all pairs of vertices $v,w\in V(G)$, let $G\cup vw$ be the graph with vertex set $V(G)$ and edge set $E(G)\cup\{vw\}$.

For each integer $k\geq1$, $k$-trees are the class of graphs defined recursively as follows. The complete graph $K_{k+1}$ is a \emph{$k$-tree}, and the graph obtained from a $k$-tree by adding a new vertex adjacent to each vertex of an existing $k$-clique is also a \emph{$k$-tree}. The \emph{treewidth} of a graph $G$ is the minimum $k$ such that $G$ is a spanning subgraph of a $k$-tree. For example, the graphs of treewidth one are the forests.  Graphs of treewidth two, called \emph{series-parallel}, are planar since in the construction of a $2$-tree, each new vertex can be drawn close to the midpoint of the edge that it is added onto.  Maximal outerplanar graphs are examples of $2$-trees.

\mySection{Some Special Plane Graphs}{DilationFree}

As illustrated in \figref{Eppstein}, we have the following characterisation of plane drawings with a segment between every pair of vertices. In this sense, these are the plane drawings with the least number of segments.

\begin{theorem}
\thmlabel{Eppstein}
In a plane drawing of a planar graph $G$, every pair of vertices of $G$ is connected by a segment if and only if at least one of the following conditions hold:
\begin{enumerate}
\item[\textup{(a)}] all the vertices of $G$ are collinear,
\item[\textup{(b)}] all the vertices of $G$, except for one, are collinear,
\item[\textup{(c)}] all the vertices of $G$, except for two vertices $v$ and $w$, are collinear, such that the line-segment $\overline{vw}$ passes through one of the collinear vertices,
\item[\textup{(d)}] all the vertices of $G$, except for two vertices $v$ and $w$,  are collinear, such that the line-segment $\overline{vw}$ does not intersect the line-segment containing $V(G)\setminus\{v,w\}$,
\item[\textup{(e)}] $G$ is the $6$-vertex octahedron graph \textup{(}say $V(G)=\{1,2,\dots,6\}$ and $E(G)=\{12,13,23,45,\linebreak 46,56,14,15,25,26,34,36\}$\textup{)} with the triangle $\{4,5,6\}$ inside the triangle $\{1,2,3\}$, and each of the triples $\{1,4,6\}$, $\{2,5,4\}$, $\{3,6,5\}$ are collinear.
\end{enumerate}
\end{theorem}

\begin{proof}
As illustrated in \figref{Eppstein}, in a plane graph that satisfies one of (a)-(e), every pair of vertices is connected by a segment. For the converse, consider a plane graph $G$ in which every pair of vertices is connected by a segment. Let $L$ be a maximum set of collinear vertices. Let $\hat{L}$ be the line containing $L$. Then $|L|\geq2$. If $|L|=2$, then $G=K_n$ for some $n\leq 4$, which is included in case (a), (b), or (d). Now suppose that $|L|\geq3$. 

Without loss of generality, $\hat{L}$ is horizontal. Let $S$ and $T$ be the sets of vertices respectively above and below $\hat{L}$. Assume $|S|\geq|T|$. 

If $|S|\leq 1$, then it is easily seen that $G$ is in case (a), (b), (c), or (d). Otherwise $|S|\geq 2$. Choose $v\in S$ to be the closest vertex to $\hat{L}$ (in terms of perpendicular distance), and choose $w\in S\setminus\{v\}$ to be the next closest vertex. This is possible since $G$ is finite. Let $p$ be the point of intersection between $\hat{L}$ and the line through $v$ and $w$. 

Suppose on the contrary that there are at least two vertices $x,y\in L$ on one side of $p$. Say $x$ is between $p$ and $y$. Then the segments $vy$ and $wx$ cross at a point closer to $\hat{L}$ than $v$. Since $G$ is plane, there is a vertex in $S$ at this point, contradicting our choice of $v$. Hence there is at most one vertex in $L$ on each side of $p$. Since $|L|\geq3$, $p$ is a vertex in $L$, and $|L|=3$. Thus there is exactly one vertex in $L$ on each side of $p$. Let these vertices be $x$ and $y$. 

Suppose on the contrary that there is a vertex $u\in S\setminus\{v,w\}$. Then $u$ is above $w$, and $u$ is not on the line containing $v,w,p$ (as otherwise $L$ is not a maximum set of collinear points). Thus the segment $uv$ crosses either $wx$ or $wy$ at a point closer to $\hat{L}$ than $w$. Since $G$ is plane, there is a vertex in $S$ at this point, contradicting our choice of $w$. Thus $|S|=2$, which implies $|T|\leq 2$. 

Now $V(G)=\{v,w,p,x,y\}\cup T$. We have $|V(G)|\geq6$, as otherwise $G$ is in case (a), (b), (c) or (d). Hence $T\ne\emptyset$. Consider a vertex $q\in T$. The segment $qv$ crosses $\hat{L}$ at some vertex in $L$. It cannot cross at $p$ (as otherwise $L$ would not be a maximum set of collinear vertices). Thus every vertex $q\in T$ is collinear with $vx$ or $vy$. Suppose there are two vertices $q_1,q_2\in T$ with $q_1$ collinear with $vx$ and $q_2$ collinear with $vy$. Then the segments $q_1y$ and $q_2x$ would cross at a point below $\hat{L}$ but not collinear with $vx$ or $vy$, which is a contradiction. Suppose there are two vertices $q_1,q_2\in T$ both collinear with $vx$; say $q_1$ is closer to $\hat{L}$ than $q_2$. Then the segments $q_1y$ and $q_2p$ would cross at a point below $\hat{L}$ but not collinear with $vx$ or $vy$, which is a contradiction. We obtain a similar contradiction if there are two vertices $q_1,q_2\in T$ both collinear with $vy$. Thus there is exactly one vertex $q\in T$. Without loss of generality, $q$ is collinear with $vx$.  Then $\{v,w,x,y,p,q\}$ induce the octahedron in case (e) where $1=q$, $2=y$, $3=w$, $4=x$, $5=p$, and $6=v$. \end{proof}

\Figure{Eppstein}{\includegraphics{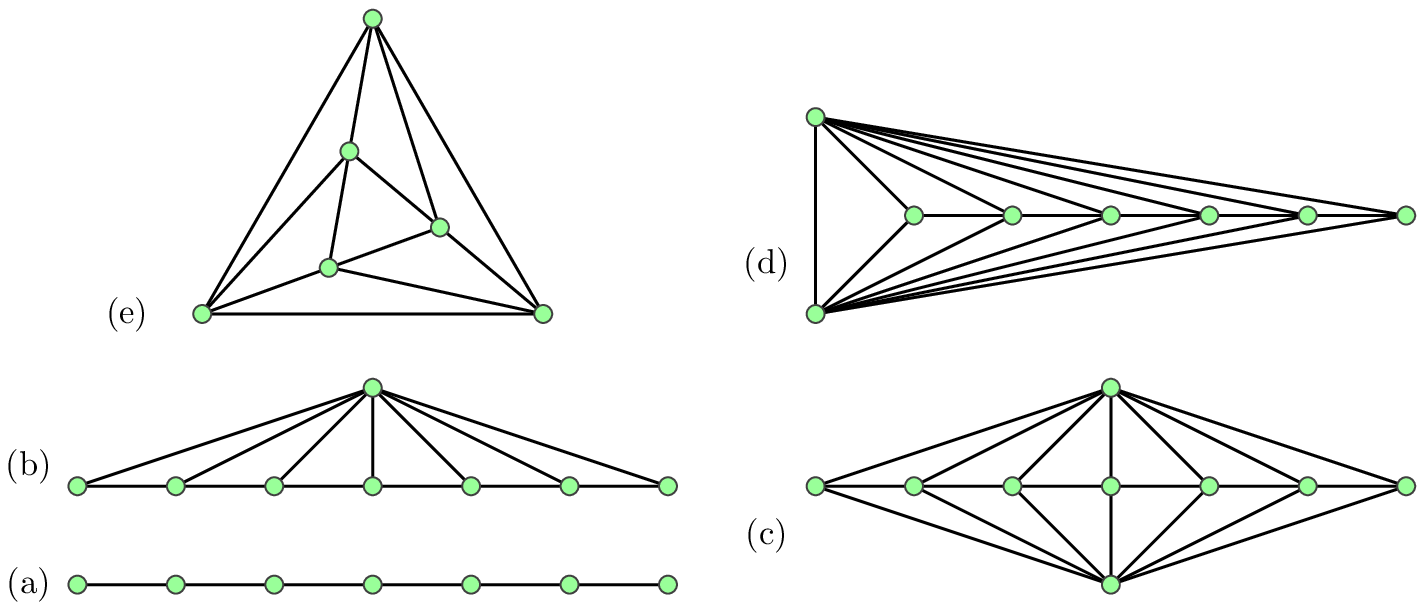}}{The plane graph drawings with a segment between every pair of  vertices.}





%
%

\mySection{Drawings on Two or Three Slopes}{TwoThree}

For drawings on two or three slopes the choice of slopes is not important.

\begin{lemma} 
\lemlabel{ThreeSlope}
A graph has a \paran{plane} drawing on three slopes if and only if it has a \paran{plane} drawing on any three slopes.
\end{lemma}

\begin{proof}
Let $D$ be a drawing of a graph $G$ on slopes $s_1,s_2,s_3$. Let $t_1,t_2,t_3$ be three given slopes. Let $T$ be a triangle with slopes $s_1,s_2,s_3$. Let $T'$ be a triangle with slopes $t_1,t_2,t_3$. It is well known that there is an affine transformation $\alpha$ to transform $T$ into $T'$. Let $D'$ be the result of applying $\alpha$ to $D$. Since parallel lines are preserved under $\alpha$, every edge in $D'$ has slope in $\{t_1,t_2,t_3\}$. Since sets of collinear points are preserved under $\alpha$, no edge passes through another vertex in $D'$. Thus $D'$ is a drawing of $G$ with slopes $t_1,t_2,t_3$. Moreover, two edges cross in $D$ if and only if they cross in $D'$. Thus $D'$ is plane whenever $D$ is plane.
\end{proof}

\begin{corollary} 
\corlabel{TwoSlope}
A graph has a \paran{plane} drawing on two slopes if and only if it has a \paran{plane} drawing on any two slopes.\qed
\end{corollary}


\citet{GT-SJC01} proved that it is \NP-complete to decide whether  a graph has a rectilinear planar drawing (that is, with vertical and horizontal edges). Thus \corref{TwoSlope} implies:

\begin{corollary}
It is \NP-complete to decide whether a graph has a plane drawing with two  slopes.\qed
\end{corollary}

Note that it is easily seen that $K_4$ has a drawing on four slopes, but does not have a drawing on slopes $\{0,\epsilon,\frac{\pi}{2},\frac{\pi}{2}+\epsilon\}$ for small enough $\epsilon$.

\mySection{Planar Graphs with Small Treewidth}{PlanarTreewidth}

\mySubSection{Trees}{Trees}

In this section we study drawings of trees with few slope and few segments. We start with the following universal lower bounds.

\begin{lemma}
\lemlabel{LowerBounds} 
The number of slopes in a drawing of a graph is at least half the maximum degree, and at least the minimum degree. The number of segments in a drawing of a graph is at least half the number of odd degree vertices.
\end{lemma}

\begin{proof} 
At most two edges incident to a vertex $v$ can have the same slope. Thus the edges incident to $v$ use at least $\half\deg(v)$ slopes. Hence the number of slopes is at least half the maximum degree. For some vertex $v$ on the convex hull, every edge incident to $v$ has a distinct slope. Thus the number of slopes is at least the minimum degree. 

If a vertex is internal on every segment then it has even degree. Thus each vertex of odd degree is an endpoint of some segment. Thus the number of vertices of odd degree is at most twice the number of segments. (The number of odd degree vertices is always even.)
\end{proof}

We now show that the lower bounds in \lemref{LowerBounds} are tight for trees. In fact, they can be simultaneously attained by the same drawing. 

\begin{theorem}
\thmlabel{Trees}
Let $T$ be a tree with maximum degree $\Delta$, and with $\eta$ vertices of odd degree.  The minimum number of segments in a drawing of $T$ is $\frac{\eta}{2}$. The minimum number of slopes in a drawing of $T$ is $\ceil{\frac{\Delta}{2}}$. Moreover, $T$ has a plane drawing with $\frac{\eta}{2}$ segments and $\ceil{\frac{\Delta}{2}}$ slopes.
\end{theorem}

\begin{proof} 
The lower bounds are \lemref{LowerBounds}. The upper bound will follow from the following hypothesis, which we prove by induction  on the number of vertices: ``Every tree $T$ with maximum degree $\Delta$ has a plane drawing with $\ceil{\frac{\Delta}{2}}$ slopes, in which every odd degree vertex is an endpoint of exactly one segment, and no even degree vertex is an endpoint of a segment.'' The hypothesis is trivially true for a single vertex. Let $x$ be a leaf of $T$ incident to the edge $xy$. Let $T'=T\setminus x$. Suppose $T'$ has maximum degree $\Delta'$.

First suppose that $y$ has even degree in $T$, as illustrated in \figref{Trees}(a). Thus $y$ has odd degree in $T'$. By induction, $T'$ has a plane drawing with $\ceil{\frac{\Delta'}{2}}\leq\ceil{\frac{\Delta}{2}}$ slopes, in which $y$ is an endpoint of exactly one segment. That segment contains some edge $e$ incident to $y$. Draw $x$ on the extension of $e$ so that there are no crossings. In the obtained drawing $D$, the number of slopes is unchanged, $x$ is an endpoint of one segment, and $y$ is not an endpoint of any segment. Thus $D$ satisfies the hypothesis.

\Figure{Trees}{\includegraphics{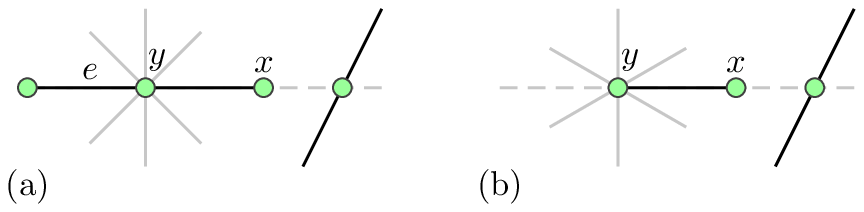}}{Adding a leaf $x$ to a drawing of a
tree: (a) $\deg(y)$ even and (b) $\deg(y)$ odd.}

Now suppose that $y$ has odd degree in $T$, as illustrated in \figref{Trees}(b). Thus $y$ has even degree in $T'$. By induction, $T'$ has a plane drawing with $\ceil{\frac{\Delta'}{2}}$ slopes, in which $y$ is not an endpoint of any segment. Thus the edges incident to $y$ use $\half\deg_{T'}(y)\leq\ceil{\frac{\Delta}{2}}-1$ slopes. If the drawing of $T'$ has any other slopes, let $s$ be one of these slopes, otherwise let $s$ be an unused slope. Add edge $xy$ to the drawing of $T'$ with slope $s$ so that there are no crossings. In the obtained drawing $D$, there is a new segment with endpoints $x$ and $y$. Since both $x$ and $y$ have odd degree in $T$, and since $x$ and $y$  were not endpoints of any segment in the drawing of $T'$, the induction hypothesis  is maintained. The number of slopes in $D$ is $\max\{\ceil{\frac{\Delta'}{2}},\half\deg_{T'}(y)+1\}\leq\ceil{\frac{\Delta}{2}}$. \end{proof}

\mySubSection{Outerplanar Graphs}{Outerplanar}

A planar graph $G$ is \emph{outerplanar} if $G$ admits a combinatorial embedding with all the vertices on the boundary of a single face.  An outerplanar graph $G$ is \emph{maximal} if $G\cup vw$ is not outerplanar for any pair of non-adjacent vertices $v,w\in V(G)$. A plane graph is \emph{outerplanar} if all the vertices are on the boundary of the outerface.   A maximal outerplanar graph has a unique outerplanar embedding.

\begin{theorem}
\thmlabel{OuterplanarSegments}
Every $n$-vertex maximal outerplanar graph $G$ has an outerplanar drawing with  at most $n$ segments.  For all $n\geq 3$, there is an $n$-vertex maximal outerplanar graph that has at least $n$ segments in any drawing.
\end{theorem}

\begin{proof} We prove the upper bound by induction on $n$ with the additional invariant that the drawing is \emph{star-shaped}.  That is, there is a point $p$ in (the interior of) some internal face of $D$, and every ray from $p$ intersects the boundary of the outerface in exactly one point. 

For $n=3$, $G$ is a triangle, and the invariant holds by taking $p$ to be any point in the internal face. Now suppose $n>3$. It is well known that $G$ has a degree-$2$ vertex $v$ whose neighbours $x$ and $y$ are adjacent, and $G'=G\setminus v$ is maximal outerplanar. By induction, $G'$ has a drawing $D'$ with at most $n-1$ segments, and there is a point $p$ in some internal face of $D'$, such that every ray from $p$ intersects the boundary of $D'$ in exactly one point. The edge $xy$ lies on the boundary of the outerface and of some internal face $F$. Without loss of generality, $xy$ is horizontal in $D'$, and $F$ is below $xy$. Since $G'$ is maximal outerplanar, $F$ is bounded by a triangle $rxy$.

For three non-collinear points $a$, $b$ and $c$ in the plane, define the \emph{wedge} $(a,b,c)$ to be the infinite region that contains the interior of the triangle $abc$, and is enclosed on two sides by the ray from $b$ through $a$ and the ray from $b$ through $c$.  By induction, $p$ is in the wedge $(y,x,r)$ or in the wedge $(x,y,r)$. By symmetry we can assume that $p$ is in $(y,x,r)$.  

Let $R$ be the region strictly above $xy$ that is contained in the wedge $(x,p,y)$. The line extending the edge $xr$ intersects $R$.  As illustrated in \figref{Outerplanar2}, place $v$ on  any point in $R$ that is on the line extending $xr$. Draw the two incident edges $vx$ and $vy$ straight.  This defines our drawing $D$ of $G$. By induction, $R\cap D'=\emptyset$. Thus $vx$ and $vy$ do not create crossings in $D$.  Every ray from $p$ that intersects $R$, intersects the boundary of $D$ in exactly one point. All other rays from $p$ intersect the same part of the boundary of $D$ as in $D'$. Since $p$ remains in some internal face, $D$ is star-shaped. By induction, $D'$ has $n-1$ segments. Since $vx$ and $rx$ are in the same segment, there is at most one segment in $D\setminus D'$. Thus $D$ is a star-shaped outerplanar drawing of $G$ with $n$ segments. This concludes the proof of the upper bound.

\Figure{Outerplanar2}{\includegraphics{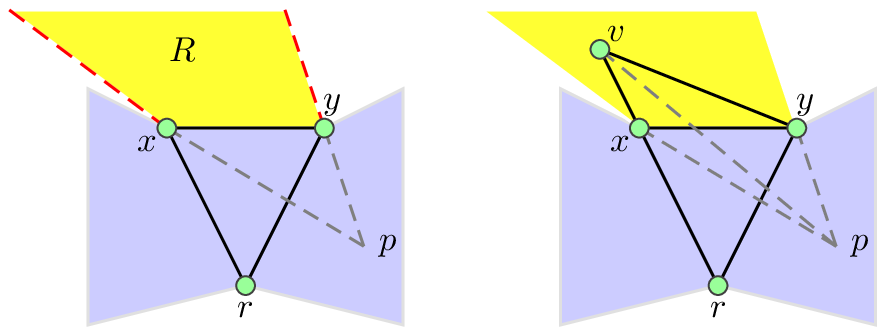}}{Construction of a
star-shaped drawing of an outerplanar graph.}

For the lower bound, let $G_n$ be the maximal outerplanar graph on $n\geq3$ vertices whose weak dual (that is, dual graph disregarding the outerface) is a path and the maximum degree of $G_n$ is at most four, as illustrated in \figref{OuterplanarLowerBound}.  

\Figure{OuterplanarLowerBound}{\includegraphics{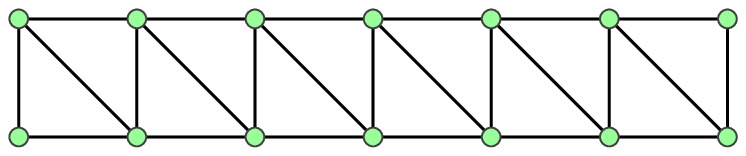}}{The
graph $G_{14}$.}

We claim that every drawing of $G_n$  has at least $n$ segments (even if crossings are allowed).  We proceed by induction on $n$. The result is trivial for $n=3$. Suppose that every drawing of $G_{n-1}$ has at least $n-1$ segments, but there exists a drawing $D$ of $G_n$ with at most $n-1$ segments.  Let $v$ be a degree-2 vertex in $G_n$ adjacent to $x$ and $y$. One of $x$ and $y$, say $x$, has degree three in $G_n$. Observe that $G_n\setminus v$ is isomorphic to $G_{n-1}$. Thus we have a drawing of $G_n$ with exactly $n-1$ segments, which contains a drawing of $G_n\setminus v$ with $n-1$ segments. Thus the edge $vx$ shares a segment with some other edge $xr$, and the edge $vy$ shares a segment with some other edge $ys$. Since $vxy$ is a triangle, $r\ne y$, $s\ne x$ and $r\ne s$. Since $x$ has degree three, $y$ is adjacent to $r$, as illustrated in \figref{MoveOuterplanar}. That accounts for all edges incident to $y$ and $x$. Thus $xy$ is a segment in $D$. 

Now construct a drawing $D'$ of $G_{n-1}$ with $x$ moved to the position of $v$ in the drawing of $G_n$.  The drawing $D$ consists of $D'$ plus the edge $xy$. Since $xy$ is a segment in $D$, $D'$ has one less segment than $D$.  Thus $D'$ is a drawing of $G_{n-1}$ with at most $n-2$ segments,  which is the desired contradiction. \end{proof}

\Figure{MoveOuterplanar}{\includegraphics{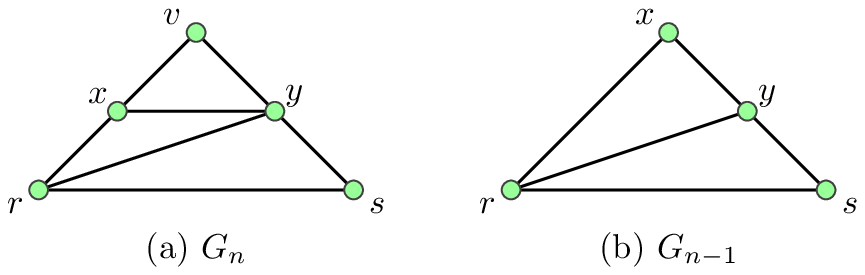}}{Construction of a
drawing of $G_{n-1}$ from a drawing of $G_n$.}

\begin{open} 
Is there a polynomial time algorithm to compute an outerplanar drawing of a given outerplanar graph with the minimum number of segments?
\end{open}

\mySubSection{$2$-Trees}{TwoTrees}

In this section we study drawings of $2$-trees with few slopes and segments. The following lower bound on the number of slopes is immediate, as illustrated in \figref{Triangle}.

\begin{observation} 
\obslabel{*}
Let $u$, $v$ and $w$ be three non-collinear vertices  in a drawing $D$ of a graph $G$. Let $d(u)$ denote the number of edges incident to $u$ that intersect the interior of the triangle $uvw$, and similarly for $v$ and $w$. Then $D$ has at least $d(u)+d(v)+d(w)+|E(G)\cap\{uv,vw,uw\}|$ slopes. \qed
\end{observation}

\Figure{Triangle}{\includegraphics{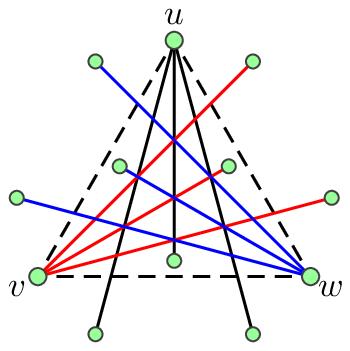}}{A triangle forces many different slopes.}

\begin{lemma}
\lemlabel{TwoTreeSegments}
Every $n$-vertex $2$-tree has a plane drawing with at most $2n-3$ segments \paran{and thus at most $2n-3$ slopes}. For all $n\geq3$, there is an $n$-vertex plane $2$-tree that has at least $2n-3$ slopes \paran{and thus at least $2n-3$ segments} in every plane drawing. 
\end{lemma}

\begin{proof} 
The upper bound follows from the F{\'a}ry-Wagner theorem  since every $2$-tree is planar and has $2n-3$ edges. Consider the $2$-tree $G_n$ with vertex set $\{v_1,v_2,\dots,v_n\}$ and edge set $\{v_1v_2,v_1v_i,v_2v_i:3\leq i\leq n\}$. Fix a plane embedding of $G_n$ with the edge $v_1v_2$ on the triangular outerface, as illustrated in \figref{TwoTree}(a). The number of slopes is at least $(n-3)+(n-3)+0+3=2n-3$ by \obsref{*}. 
\end{proof}

\Figure{TwoTree}{\includegraphics{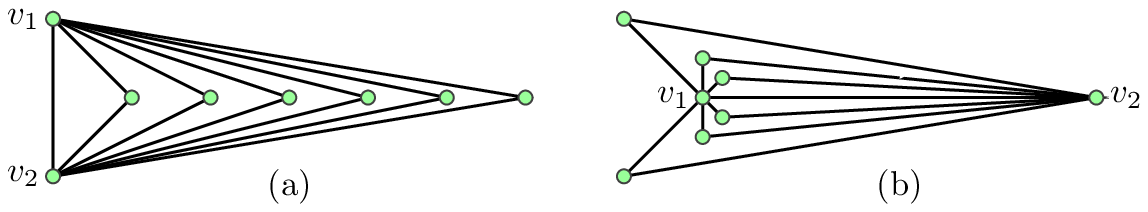}}{The graph $G_8$ in \lemref{TwoTreeSegments}.}

In \lemref{TwoTreeSegments} the embedding is fixed. A better bound can be obtained if we do not fix the embedding. For example, the graph $G_n$ from \lemref{TwoTreeSegments} has a plane drawing with $\frac{3n}{2}-2$ segments, as illustrated in \figref{TwoTree}(b). 

\begin{theorem}
\thmlabel{TwoTreeSegmentsImproved}
Every $n$-vertex $2$-tree $G$ has a plane drawing with at most $\frac{3}{2}n$ segments
\end{theorem}

The key idea in the proof of \thmref{TwoTreeSegmentsImproved} is to position a set of vertices at each step, rather than a single vertex. The next lemma says how to partition a $2$-tree appropriately. It has subsequently been generalised for $k$-trees by \citet{DujWoo-GD05}.

\begin{lemma}
\lemlabel{TwoTreePartition}
Let $G$ be a $2$-tree. Then for some $k\geq1$, $V(G)$ can be partitioned $(S_0,S_1,S_2,\dots,S_k)$ such that:
\begin{enumerate}
\item[\textup{(a)}] for $1\leq i\leq k$, the subgraph $G_i=G[\bigcup_{j=0}^iS_j]$ is a $2$-tree,
\item[\textup{(b)}] $S_0$ consists of two adjacent vertices, 
\item[\textup{(c)}] for $1\leq i\leq k$, $S_i$ is an independent set of $G$,
\item[\textup{(d)}] for $1\leq i\leq k$, each vertex in $S_i$ has exactly two neighbours in $G_{i-1}$, and they are adjacent,
\item[\textup{(e)}] for $2\leq i\leq k$, the vertices in $S_i$ have a common neighbour $v$ in 
$G_{i-1}$, and $v$ has degree two in $G_{i-1}$.
\end{enumerate}
\end{lemma}

\begin{proof}
We proceed by induction on $|V(G)|$. By definition, $|V(G)|\geq 3$. First suppose that $|V(G)|=3$. Let $V(G)=\{u,v,w\}$. Then $G=K_3$, and $(\{u,v\},\{w\})$ is the desired partition of $G$. Now suppose that $|V(G)|>3$. Let $L$ be the set of vertices of degree two in $G$. Then $L$ is a nonempty independent set, the neighbours of each vertex in $L$ are adjacent, and $G\setminus L$ is a $2$-tree. If $G\setminus L$ is a single edge $vw$, then $(\{v,w\},L)$ is the desired partition of $G$. Otherwise, $G\setminus L$ has a vertex $v$ of degree two. Let $S$ be the set of neighbours of $v$ in $L$. Now $S\ne\emptyset$, as otherwise $v\in L$. By induction, there is a partition $(S_0,S_1,S_2,\dots,S_k)$ of $V(G\setminus S)$ that satisfies the lemma. It is easily verified that $(S_0,S_1,S_2,\dots,S_k,S)$ is the desired partition of $G$. 
\end{proof}


\begin{proof}[Proof of \thmref{TwoTreeSegmentsImproved}]
Let $(S_0,S_1,S_2,\dots,S_k)$ be the partition of $V(G)$ from \lemref{TwoTreePartition}. First suppose that $k=1$. By \lemref{TwoTreePartition}(b) and (d), $S_0=\{v,w\}$ and $S_1$ is an independent set of vertices, each adjacent to both $v$ and $w$. Let $S_1=\{a_1,a_2,\dots,a_p\}\cup\{b_1,b_2,\dots,b_q\}$, where $q\leq p\leq q+1$. As illustrated in \figref{TwoTreeSegmentsImproved}(a), $G$ can be drawn such that $a_iv$ and $b_iv$ form a single segment, for all $1\leq i\leq q$.
The number of segments is at most $1+|S_1|+\ceil{\half|S_1|}\leq\half(3n-3)$.

Now suppose that $k\geq2$. By \lemref{TwoTreePartition}(a), $G_{k-1}$ is a $2$-tree. Thus by induction, $G_{k-1}$ has a plane drawing with at most $\frac{3}{2}(n-|S_k|)$ segments. By \lemref{TwoTreePartition}(d) and (e), the vertices in $S_k$ have degree two in $G$, and have a common neighbour $v$ in $G_{k-1}$ with degree two in $G_{k-1}$. Let $u$ and $w$ be the neighbours of $v$ in $G_{k-1}$. Then the neighbourhood of each vertex in $S_k$ is either $\{v,u\}$ or $\{v,w\}$. Let $S_k^u$ and $S_k^w$ be the sets of vertices in $S_k$ whose neighbourhood respectively is $\{v,u\}$ and $\{v,w\}$. Without loss of generality, $|S_k^u|\geq|S_k^w|$. Let $S_k^w=\{b_1,\dots,b_p\}$. For the time being assume that $|S_k^u|-p$ is even. Let $r=\half(|S_k^u|-p)$. Thus $r$ is a nonnegative integer, and $S_k^u$ can be partitioned 
\begin{equation*}
S_k^u=\{a_1,\dots,a_p\}\cup\{c_1,\dots,c_r\}\cup\{d_1,\dots,d_r\}\enspace.
\end{equation*}

As illustrated in \figref{TwoTreeSegmentsImproved}(b), $G$ can be drawn such that $a_iv$ and $b_iv$ form a single segment for all $1\leq i\leq p$, and
$c_iv$ and $d_iv$ form a single segment for all $1\leq i\leq r$. 
Clearly the vertices can be placed to avoid crossings with the existing drawing of $G_{k-1}$. In particular, vertices $\{b_1,\dots,b_p,d_1,\dots,d_r\}$ are drawn inside the triangle $(u,v,w)$. The number of new segments in the drawing is $3p+3r=\frac{3}{2}|S_k|$. 

In the case that $|S_k^u|-p$ is odd, a vertex $x$ from $S_k^u$ can be drawn so that $xv$ and $xw$ form a single segment; then apply the above algorithm to $S_k\setminus\{x\}$. The number of new segments is then $3p+3r+1$, where $|S_k|=2p+2r+1$. It follows that the number of new segments is at most $\half(3|S_k|-1)$. 

In both cases, the total number of segments is at most $\frac{3}{2}(n-|S_k|)+\frac{3}{2}|S_k|=\frac{3}{2}n$. 
\end{proof}

\Figure{TwoTreeSegmentsImproved}{\includegraphics{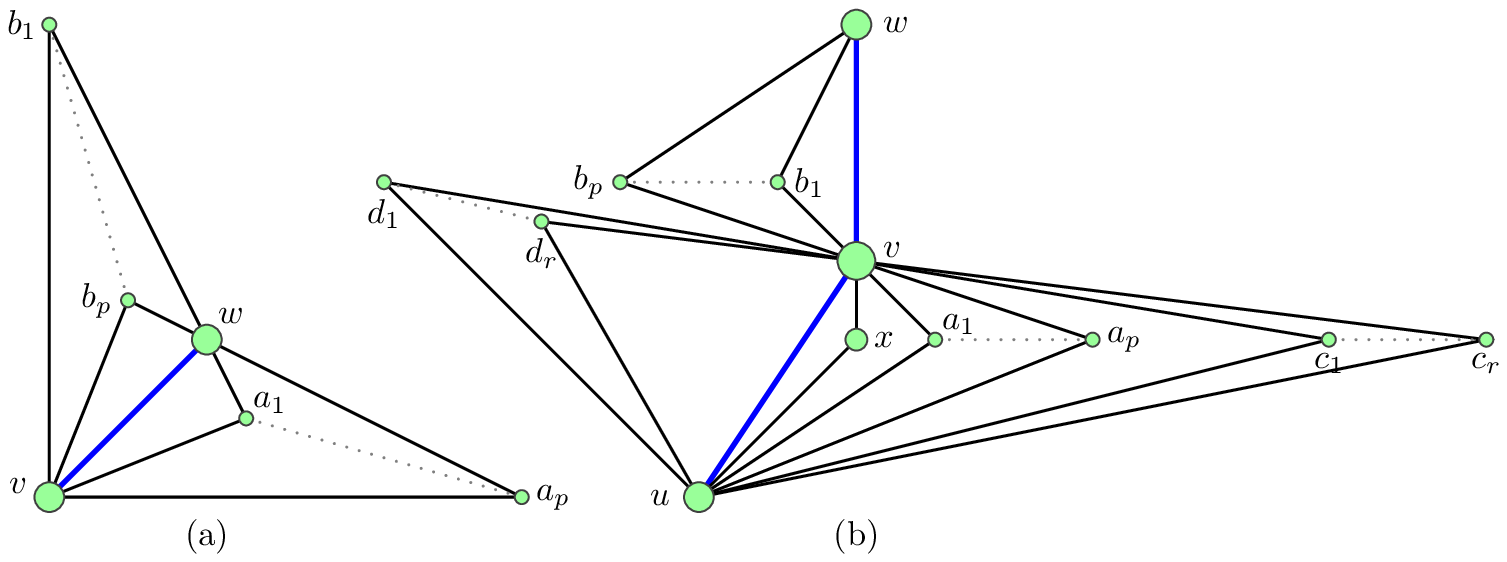}}{(a) drawing $S_1$, (b) drawing $S_k$.}

\mySubSection{Planar $3$-Trees}{ThreeTrees}

We now turn our attention to drawings of planar $3$-trees. 

%
%
%

\begin{theorem}
\thmlabel{PlanarThreeTreeSegments}
Every $n$-vertex plane $3$-tree has a plane drawing with at most $2n-2$ segments \paran{and thus at most $2n-2$ slopes}. For all $n\geq4$, there is an $n$-vertex plane $3$-tree with at least $2n-2$ slopes \paran{and thus at least $2n-2$ segments} in every drawing.
\end{theorem}

\begin{proof} We prove the upper bound by induction on $n$ with the hypothesis that  ``every plane $3$-tree with $n\geq4$ vertices has a plane  drawing with at most $2n-2$ segments, such that for every internal face $F$ there is an edge $e$ incident to exactly one vertex of $F$, and the extension of $e$ intersects the interior of $F$.'' The base case is trivial since $K_4$ is the only $3$-tree on four vertices, and any plane drawing of $K_4$ satisfies the hypothesis.

Suppose that the claim holds for plane $3$-trees on $n-1$ vertices. Let $G$ be a plane  $3$-tree on $n$ vertices. Every $k$-tree on at least $k+2$ vertices has two non-adjacent  simplicial vertices of degree exactly $k$ \citep{Dirac61}. In particular, $G$ has two  non-adjacent simplicial degree-$3$ vertices, one of which, say $v$, is not on the outerface. Thus $G$ can be obtained from $G\setminus v$ by adding $v$ inside some internal face $(p,q,r)$ of $G\setminus v$, adjacent to $p$, $q$ and $r$\footnote[5]{Note that this implies that the planar $3$-trees are precisely those graphs that are produced by the LEDA `random' maximal planar graph generator. This algorithm, starting from $K_3$, repeatedly adds a new vertex adjacent to the three vertices of a randomly selected internal face.}.  By induction, $G\setminus v$ has a drawing with $2n-4$ segments in which there is an edge $e$ incident to exactly one of $\{p,q,r\}$, and the extension of $e$ intersects the interior of the face. Position $v$ in the interior of the face anywhere on the extension of $e$, and draw segments from $v$ to each of $p$, $q$ and $r$. We obtain a plane drawing of $G$ with $2n-2$ segments. The extension of $vp$ intersects the interior of $(v,q,r)$; the extension of $vq$ intersects the interior of $(v,p,r)$;  and the extension of $vr$ intersects the interior of $(v,p,q)$. All other faces of $G$ are faces of $G\setminus v$. Thus the inductive hypothesis holds for $G$, and the proof of the upper bound is complete.

For each $n\geq4$ we now provide a family $\mathcal{G}_n$ of $n$-vertex plane $3$-trees, each of which require at least $2n-2$ segments in any drawing. Let $\mathcal{G}_4=\{K_4\}$.  Obviously every plane drawing of $K_4$ has six segments. For all $n\geq5$,  let $\mathcal{G}_n$ be the family of plane $3$-trees $G$ obtained from some plane $3$-tree $H\in\mathcal{G}_{n-1}$ by adding a new vertex $v$ in the outerface of $H$ adjacent to each of the three vertices of the outerface.  Any drawing of $G$ contains a drawing of $H$, which contributes at least $2n-4$ segments by induction. In addition, the two edges incident to $v$ on the triangular outerface of $G$ are each in their own segment. Thus $G$ has at least $2n-2$ segments. \end{proof}

\mySection{$3$-Connected Plane Graphs}{ThreeConnected}

The following is the main result of this section.

\begin{theorem}
\thmlabel{ThreeConnectedSegmentsSlopes}
Every $3$-connected plane graph with  $n$ vertices has a plane  drawing with at most $\frac{5}{2}n-3$ segments and at most $2n-10$ slopes.
\end{theorem}

The proof of \thmref{ThreeConnectedSegmentsSlopes} is based on the canonical ordering of \citet{Kant96}, which is a generalisation of a similar structure for plane triangulations introduced by \citet{dFPP90}. Let $G$ be a $3$-connected plane graph.  \citet{Kant96} proved that $G$ has a canonical ordering defined as follows. Let $\sigma=(V_1,V_2,\dots, V_K)$ be an ordered partition of $V(G)$. That is, $V_1\cup V_2\cup\dots\cup V_K = V(G)$ and $V_i\cap V_j =\emptyset$ for all $i\ne j$. Define $G_i$ to be the plane subgraph of $G$ induced by $V_1\cup V_2\cup\dots\cup V_i$. Let $C_i$ be the subgraph of $G$ induced by the edges on the boundary of the outerface of $G_i$. As illustrated in \figref{CanonicalOrdering}, $\sigma$ is a \emph{canonical ordering} of $G$ (also called a \emph{canonical decomposition}) if:

\begin{itemize}
\item $V_1=\{v_1,v_2\}$, where $v_1$ and $v_2$ lie on the outerface and $v_1v_2\in E(G)$.
\item $V_K=\{v_n\}$, where $v_n$ lies on the outerface, $v_1 v_n\in E(G)$, and $v_n\ne v_2$. 
\item Each $C_i \ ( i>1)$ is a cycle containing $v_1v_2$.
\item Each $G_i$ is biconnected and internally $3$-connected; that is, removing any two interior vertices of $G_i$ does not disconnect it.
\item For each $i\in\{2,3,\dots,K-1\}$, one of the following conditions holds: 
\begin{enumerate}
\item $V_i=\{v_i\}$ where $v_i$ is a vertex of $C_i$ with at least three neighbours in $C_{i-1}$, and $v_i$ has at least one neighbour in $G\setminus G_i$.
\item $V_i=(s_1,s_2,\dots,s_\ell,v_i)$, $\ell\geq 0$, is a path in $C_i$, where each vertex in $V_i$ has at least one neighbour in $G\setminus G_i$. Furthermore, the first and the last vertex in $V_i$ have one neighbour in $C_{i-1}$, and these are the only two edges between $V_i$ and $G_{i-1}$. 
\end{enumerate}
\end{itemize}

\Figure{CanonicalOrdering}{\includegraphics{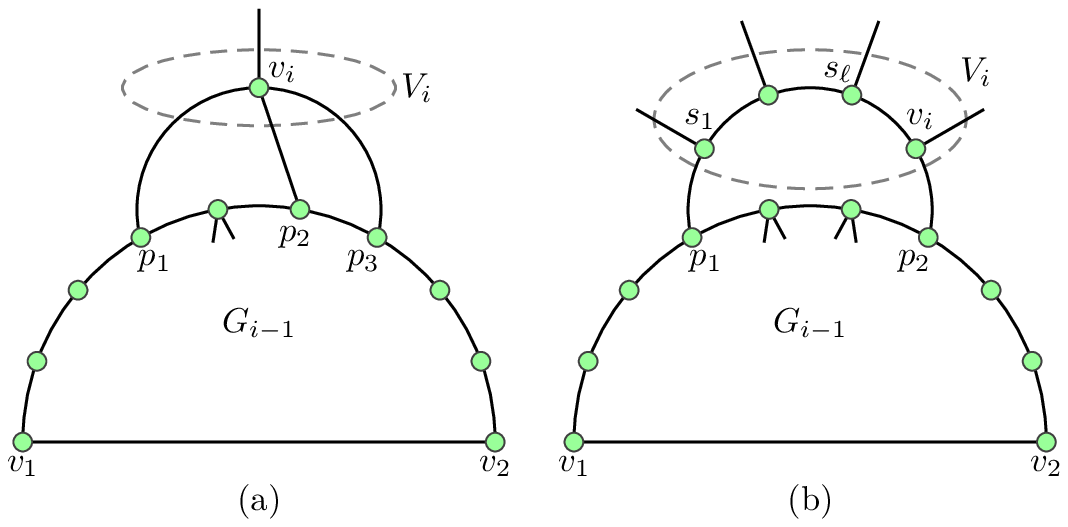}}{The canonical ordering of a $3$-connected plane graph.}

The vertex $v_i$ is called the  \emph{representative} vertex of the set $V_i$, $2\leq i \leq K$.  The vertices $\{s_1,s_2, \dots, s_\ell\}\subseteq V_i$ are called \emph{division} vertices. Let $S\subset V(G)$ be the set of all division vertices. A vertex $u$ is a \emph{successor} of a vertex $w\in V_i$ if $uw$ is an edge and $u\in G\setminus G_i$.  A vertex $u$ is a \emph{predecessor} of a vertex $w\in V_i$ if $uw$ is an edge and $u\in V_j$ for some $j<i$.  We also say that $u$ is a predecessor of $V_i$. Let $\pred{V_i}=(p_1,p_2,\dots, p_q)$ denote the set of predecessors of $V_i$ ordered by the path from $v_1$ to $v_2$ in $C_{i-1}\setminus v_1v_2$. Vertex $p_1$ and $p_q$ are the \emph{left} and \emph{right predecessors} of $V_i$ respectively, and vertices $p_2,p_3,\dots p_{q-1}$ are called \emph{middle predecessors} of $V_i$.

\begin{theorem}
\thmlabel{ThreeConnectedS}
Let $G$ be an $n$-vertex $m$-edge plane $3$-connected graph with a canonical
ordering $\sigma$. Define $S$ as above \paran{with respect to $\sigma$}. Then $G$ has a plane drawing $D$ 
with at most
\begin{equation*}
m-\max{\{\ceil{\frac{n}{2}}-|S|-3, \ |S|\}}
\end{equation*}
segments, and at most
\begin{equation*}
m-\max{\{n-|S| -4, \ |S|\}}
\end{equation*}
slopes.
\end{theorem}

\begin{proof} We first define $D$ and then determine the upper bounds on the number of segments and slopes in $D$. For every vertex $v$, let $X(v)$ and $Y(v)$ denote the $x$ and $y$ coordinates of $v$, respectively. If a  vertex $v$ has a neighbour $w$, such that $X(w)<X(v)$ and $Y(w)<Y(v)$, then we say $vw$ is a \emph{left edge} of $v$. Similarly, if $v$ has a neighbour $w$, such that $X(w)>X(v)$ and $Y(w)<Y(v)$, then we say $vw$ is a \emph{right edge} of $v$. If $vw$ is an edge such that $X(v)=X(w)$ and $Y(v)<Y(w)$, than we say $vw$ is a \emph{vertical edge above} $v$ and \emph{below} $w$.

We define $D$ inductively on $\sigma=(V_1,V_2,\dots, V_K)$ as follows. Let $D_i$ denote a drawing of $G_i$. A vertex $v$ is a \emph{peak in $D_i$}, if each neighbour $w$ of $v$ has $Y(w)\leq Y(v)$ in $D_i$. We say that a point $p$ in the plane is \emph{visible in $D_i$} from  vertex $v\in D_i$, if the segment $\overline{pv}$ does not intersect $D_i$ except at $v$. At the $i^\text{th}$ induction step, $2\leq i \leq K$, $D_i$ will satisfy the following invariants:

\begin{description}

\item[Invariant 1:] 
$C_i\setminus v_1 v_2$ is \emph{strictly X-monotone}; that is, the path from $v_1$ to $v_2$ in $C_i\setminus v_1 v_2$ has increasing X-coordinates.

\item[Invariant 2:] Every peak in $D_i$, $i<K$, has a successor.

\item[Invariant 3:] Every representative vertex $v_j\in V_j$, $2\leq j\leq i$ has a left and a right edge. Moreover, if $|\pred{V_j}|\geq 3$ then there is a vertical edge below $v_j$.

\item[Invariant 4:] $D_i$ has no edge crossings.
\end{description}

For the base case $i=2$, position the vertices $v_1$, $v_2$ and $v_3$ at the corners of an equilateral triangle so that $X(v_1)<X(v_3)<X(v_2)$ and $Y(v_1)<Y(v_2)<Y(v_3)$.  Draw the division vertices of $V_2$ on the segment $v_1v_3$. This drawing of $D_2$ satisfies all four invariants. Now suppose that we have a drawing of $D_{i-1}$ that satisfies the invariants. There are two cases to consider  in the construction of $D_i$, corresponding to the two cases in the definition of the canonical ordering.

\medskip\noindent\textbf{Case 1.} $|\pred{V_i}|\geq 3$: If $v_i$ has a middle predecessor $v_j$ with $|\pred{V_j}|\geq 3$, let $w=v_j$. Otherwise let $w$ be any middle predecessor of $v_i$. Let $L$ be the open ray $\{(X(w),y):y>Y(w)\}$. By invariant~1 for $D_{i-1}$, there is a point in $L$ that is visible in $D_{i-1}$ from every predecessor of $v_i$.  Represent $v_i$ by such a point, and draw segments between $v_i$ and each of  its predecessors. That the resulting drawing $D_i$ satisfies the four invariants can be immediately verified.

\medskip\noindent\textbf{Case 2.} $|\pred{V_i}|=2$:  Suppose that $\pred{V_i}=\{w,u\}$, where $w$ and $u$ are the left and the right predecessors of $V_i$, respectively. Suppose $Y(w)\geq Y(u)$. (The  other case is symmetric.) Let $P$ be the path between  $w$ and $u$ on $C_{i-1}\setminus  v_1v_2$.  As illustrated in \figref{CaseTwo}, let $A_i$ be the region   $\{(x,y):y>Y(w)\text{ and }X(w)\leq x\leq X(u) \}$.

\Figure{CaseTwo}{\includegraphics{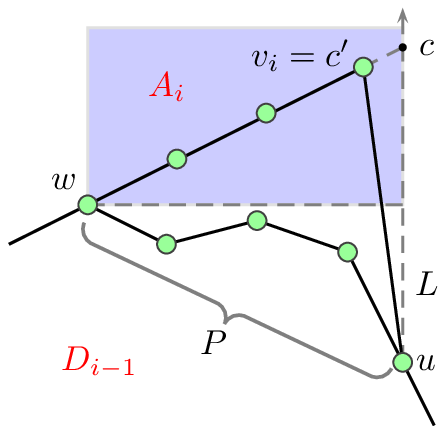}}{Illustration for Case 2.}

Assume, for the sake of contradiction, that $D_{i-1}\cap A_i\ne\emptyset$. By the monotonicity of $D_{i-1}$, $P\cap A_i \ne\emptyset$. Let $p\in P\cap A_i$. Since $Y(p)>Y(w)\geq Y(u)$, $P$ is $X$-monotone  and thus has a vertex between $w$ and $u$ that is a peak. By the definition of the canonical ordering $\sigma$, the addition of $V_i$  creates a face of $G$, since $V_i$ is added in the outerface of $G_{i-1}$. Therefore, each vertex between $w$ and $u$ on $P$ has no successor, and is thus not a peak in $D_{i-1}$ by invariant~2, which is the desired contradiction. Therefore $D_{i-1}\cap A_i=\emptyset$. 

Let $L$ be the open ray $\{(X(u),y):y>Y(u)\}$. If $w\not\in S$, then by invariant~3, $w$ has a left and a right edge in $D_{i-1}$. Let $c$ be the point of intersection between $L$ and  the line extending the left edge at $w$. If $w\in S$, then let $c$ be any point in $A_i$ on $L$. By invariant~1, there is a point $c'\not\in\{c,w\}$ on $\overline{wc}$ such that $c'$ is visible in $D_{i-1}$ from $u$. Represent $v_i$ by $c'$, and draw two segments $\overline{v_i u}$ and $\overline{v_i w}$. These two segments do not intersect any part of $D_{i-1}$ (and neither is horizontal). Represent any division vertices in $V_i$  by arbitrary  points on the open segment $\overline{w v_i}\cap A_i$. Therefore, in the resulting drawing $D_i$, there are no crossings and the remaining three invariants are maintained.  

\medskip This completes the construction of $D$. The following claim will be used to bound the number of segments and slopes in $D$. It basically says that a division vertex (and $v_2$) can be the higher predecessor for at most one set $V_i$ with $|\pred{V_i}|=2$.

\begin{claim}
Let $V_i, V_j\in \sigma$ with $i<j$ and $|\pred{V_i}|=|\pred{V_j}|=2$. Let $w_i$ be the higher of the two predecessors of $V_i$ in $D_{i-1}$, and let $w_j$ be the higher of the two predecessors of $V_j$ in $D_{j-1}$. If $w_i\in S$ or $w_i=v_2$, then $w_i\ne w_j$. 
\end{claim}

\begin{proof} Suppose that $w_i\in V_k$, $k<i$. First assume that $w_i\in S$. Then each division vertex lies on some non-horizontal segment and it is not an endpoint of that segment. Thus $w_i$ is not a peak in $D_k$, and therefore it is not a peak in every $D_\ell$, $\ell\geq k$. For all $\epsilon>0$, let
\begin{align*}
A'_\epsilon&=\{(x,y)\,:\, y> Y(w_i),\  X(w_i)-\epsilon\leq x< X(w_i)\},
\text{ and }\\
A''_\epsilon&=\{(x,y)\,:\, y> Y(w_i),\  X(w_i)<x\leq X(w_i)+\epsilon \}\enspace.
\end{align*}

Then for all small enough $\epsilon$, either $A'_\epsilon\cap D_k\ne \emptyset$ or $A''_\epsilon\cap D_k\ne \emptyset$. Without loss of generality, $A'_\epsilon\cap D_k= \emptyset$ and  $A''_\epsilon\cap D_k\ne \emptyset$. Then at iteration $i>k$, the region $A_i$, as defined in  Case~2 of the construction of $D_i$, contains $A'_\epsilon$ for all small enough $\epsilon$. Thus, $A'_\epsilon\cap D_{i} \ne\emptyset$ for all small enough $\epsilon$. Since $j\geq i+1$, $A'_\epsilon\cap D_{j-1}\ne \emptyset$ or $A''_\epsilon\cap D_{j-1}\ne \emptyset$ for all small enough $\epsilon$. Therefore, $w_i\ne w_j$ (since $V_j$ is drawn by Case~2 of the construction of $D_j$, where it is known that $A_j\cap D_{j-1}=\emptyset$). The case $w_i= v_2$ is the same, since the region $A''_\epsilon\cap D_i=\emptyset$, for every $\epsilon$ and every $1\leq i\leq K$, so only region $A'_\epsilon$ is used, and thus the above argument applies.  \end{proof}

For the purpose of counting the number of segments and
slopes in $D$ assume that we draw edge $v_1v_2$ at iteration step $i=1$ and $G_2\setminus v_1v_2$ at iteration $i=2$. In every iteration $i$ of the construction, $2\leq i\leq K$, at most $|\pred{V_i}|$ new segments and slopes are created. We call an iteration $i$ of the construction \emph{segment-heavy} if the difference between the number of segments in $D_i$ and $D_{i-1}$ is exactly $|\pred{V_i}|$, and  \emph{slope-heavy} if the difference between the number of slopes in $D_i$ and $D_{i-1}$ is exactly $|\pred{V_i}|$.  Let $h_s$ and $h_\ell$ denote the total number of segment-heavy and slope-heavy iterations, respectively. Then $D$ uses at most 
\begin{equation}
\eqnlabel{seg}
1+\sum_{i=2}^K (|\pred{V_i}|-1) +h_s
\end{equation}  
segments,  and at most 
\begin{equation}
\eqnlabel{slo}
1+\sum_{i=2}^K (|\pred{V_i}|-1) +h_\ell
\end{equation}
 slopes. 

We first express $\sum_{i=2}^K |\pred{V_i}|$ in terms of $m$ and $|S|$, and then establish an upper bound on $h_s$ and $h_\ell$.  For $i\geq 2$, let $E_i$ denote the set of edges of $G_i$ with at least one endpoint in $V_i$, and let $\ell_i$ denote the number of division vertices in $V_i$. Then $m=1+\sum_{i=2}^K |E_i| =  1+\sum_{i= 2}^K(\ell_i +|\pred{V_i}|)=  1+|S| + \sum_{i=2}^K |\pred{V_i}|$. Thus $\sum_{i=2}^K |\pred{V_i}|=m-|S|-1$. Since the trivial upper bound for $h_s$ and $h_\ell$ is $K-1$, and by \eqnref{seg} and \eqnref{slo}, we have that $D$ uses at most  $1+\sum_{i=2}^K |\pred{V_i}|= 1+ m-|S|-1 =m-|S|$ segments and slopes.

We now prove a tighter bound on $h_s$.  Let $R$ denote the set of representative vertices of   segment-heavy steps $i$ with $|\pred{V_i}|\geq 3$. Consider a step $i$ such that $|\pred{V_i}|\geq 3$.  If $v_i$ has at least one  predecessor $v_j$ with $|\pred{V_j}|\geq 3$, then  $v_i$ is drawn on the line that extends the vertical edge below $v_j$, and thus step $i$ introduces at most $|\pred{V_i}|-1$ new segments and is not segment-heavy.    Therefore, step $i$ is segment-heavy only if no middle  predecessor $w$ of $v_i$ is in $R$. Thus for each  segment-heavy step $i$ with $|\pred{V_i}|\geq 3$, there is  a unique vertex $w\not\in R$. In other words, for each  vertex in $R$, there is a unique vertex in $V(G)\setminus  R$. Thus $|R|\leq\floor{\frac{n}{2}}$. Since the number of  segment-heavy steps $i$ with $|\pred{V_i}|\geq 3$ is equal to $|R|$, there is at most $\floor{\frac{n}{2}}$ such steps.
 

The remaining steps, those with $|\pred{V_i}|=2$, introduce $|\pred{V_i}|$ segments only if the higher of the two predecessors  of $V_i$ is in $S$ or is $v_2$. (It cannot be $v_1$, since $Y(v_1)<Y(v)$  for every vertex $v\ne v_1$.)\  By the above claim, there may be at most $|S|+1$ such segment-heavy steps. Therefore, $h_s\leq \floor{\frac{n}{2}}+|S|+1$. By \eqnref{seg} and since $K=n-1-|S|$, $D$ has at most $m-\ceil{\frac{n}{2}}+ |S| +3$ segments.

Finally, we bound $h_\ell$. There may be at most one slope-heavy step $i$ with $|\pred{v_i}|\geq 3$, since there is a vertical edge below every such vertex $v_i$ by invariant~3. As in the above case for segments, there may be at most $|S|+1$ slope-heavy steps $i$ with $|\pred{v_i}|=2$. Therefore, $h_\ell\leq |S|+2$. By \eqnref{slo} and since $K= n-1-|S|$, we have that $D$ has at most $m-n +|S|+4$ slopes. \end{proof}

\begin{proof}[Proof of \thmref{ThreeConnectedSegmentsSlopes}]  Whenever a set $V_i$ is added to $G_{i-1}$, at least $|V_i|-1$ edges that are not in $G$ can be added so that the resulting graph is planar. Thus $|S|=\sum_i(|V_i|-1)\leq3n-6-m$. Hence \thmref{ThreeConnectedS} implies that $G$ has a plane drawing  with at most $m-\frac{n}{2}+|S|+3\leq\frac{5}{2}n-3$ segments, and at most $m-n+|S|-4\leq2n-10$ slopes. \end{proof}


We now prove that the bound on the number of segments in \thmref{ThreeConnectedSegmentsSlopes} is tight.

\begin{lemma}
\lemlabel{TriangulationsSegments}
For all $n\equiv0\pmod{3}$, there is an $n$-vertex planar triangulation  with maximum degree six that has at least $2n-6$ segments in every plane drawing, regardless of the choice of outerface.
\end{lemma}

\begin{proof} Consider the planar triangulation $G_k$ with vertex set $\{x_i,y_i,z_i:1\leq i\leq k\}$ and edge set $\{x_iy_i,y_iz_i,z_ix_i:1\leq i\leq k\}\cup\{x_ix_{i+1},y_iy_{i+1},z_iz_{i+1}:1\leq i\leq k-1\}\cup \{x_iy_{i+1},y_iz_{i+1},z_ix_{i+1}:1\leq i\leq k-1\}$.  $G_k$ has $n=3k$ vertices. $G_k$ is the famous `nested-triangles' graph. We say $\{(x_i,y_i,z_i):1\leq i\leq k\}$ are \emph{the triangles} of $G_k$.  This graph has a natural plane embedding with the triangle $x_iy_iz_i$ nested inside the triangle $(x_{i+1},y_{i+1},z_{i+1})$ for all $1\leq i\leq k-1$, as illustrated in \figref{NestedTriangles}. 

\Figure{NestedTriangles}{\includegraphics{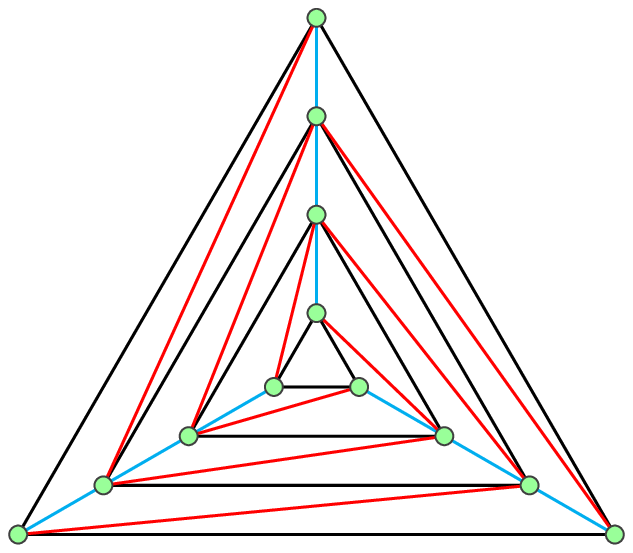}}{The graph $G_4$ in \lemref{TriangulationsSegments}.}

We first prove that if $(x_k,y_k,z_k)$ is the outerface then $G_k$ has at least $6k$ segments in any plane drawing. First observe that no two edges in the triangles can share a segment. Thus they contribute $3k$ segments. 

We claim that the six edges between triangles  $(x_i,y_i,z_i)$ and $(x_{i+1},y_{i+1},z_{i+1})$  contribute a further three segments. Consider the two edges $x_ix_{i+1}$ and $z_ix_{i+1}$ incident on $x_{i+1}$.  We will show that at least one of them contributes a new segment.  Let $R_x$ be the region bounded by the lines containing $x_iy_i$ and $x_iz_i$ that shares only $x_i$ with triangle $(x_i,y_i,z_i)$.  Similarly, let $R_z$ be the region bounded by the lines containing $x_iz_i$ and $y_iz_i$ that shares only $z_i$ with the same triangle.  We note that these two regions are disjoint.
Furthermore, if edge $x_ix_{i+1}$ belongs to a segment including edges contained in triangle $(x_i,y_i,z_i)$, then $x_{i+1}$ lies in region $R_x$.  Similarly, if $z_ix_{i+1}$ belongs to a segment including edges
contained in triangle $(x_i,y_i,z_i)$, then $x_{i+1}$ lies in region $R_z$.  Both cases cannot be true simultaneously so either edge $x_ix_{i+1}$ or edge $z_ix_{i+1}$ contributes a new segment to the drawing.  Symmetric arguments apply to the edges incident on $y_{i+1}$ and $z_{i+1}$ so the edges between triangles contribute at least three segments.

Thus in total we have at least $3k+3(k-1)=2n-3$ segments.  Now suppose that some face, other than $(x_k,y_k,z_k)$, is the outerface. Thus the triangles are split into two nested sets. Say there are $p$ triangles in one set and $q$ in the other. By the above argument, any drawing has at least $(2p-3)+(2q-3)=2n-6$ segments. \end{proof}


\lemref{TriangulationsSegments} gives a tight lower bound of $2n-c$ on the number of segments in drawings of maximal planar graphs. However, there are plane drawings with as little as \Oh{\sqrt{n}} segments, as illustrated in \figref{FewSegments}. Note that for planar graphs without degree two vertices, if there are $k$ segments in some drawing, then the corresponding arrangement has at most $\binom{k}{2}$ vertices. Thus $n\leq\binom{k}{2}$ and $k>\sqrt{2n}$.

\Figure{FewSegments}{\includegraphics{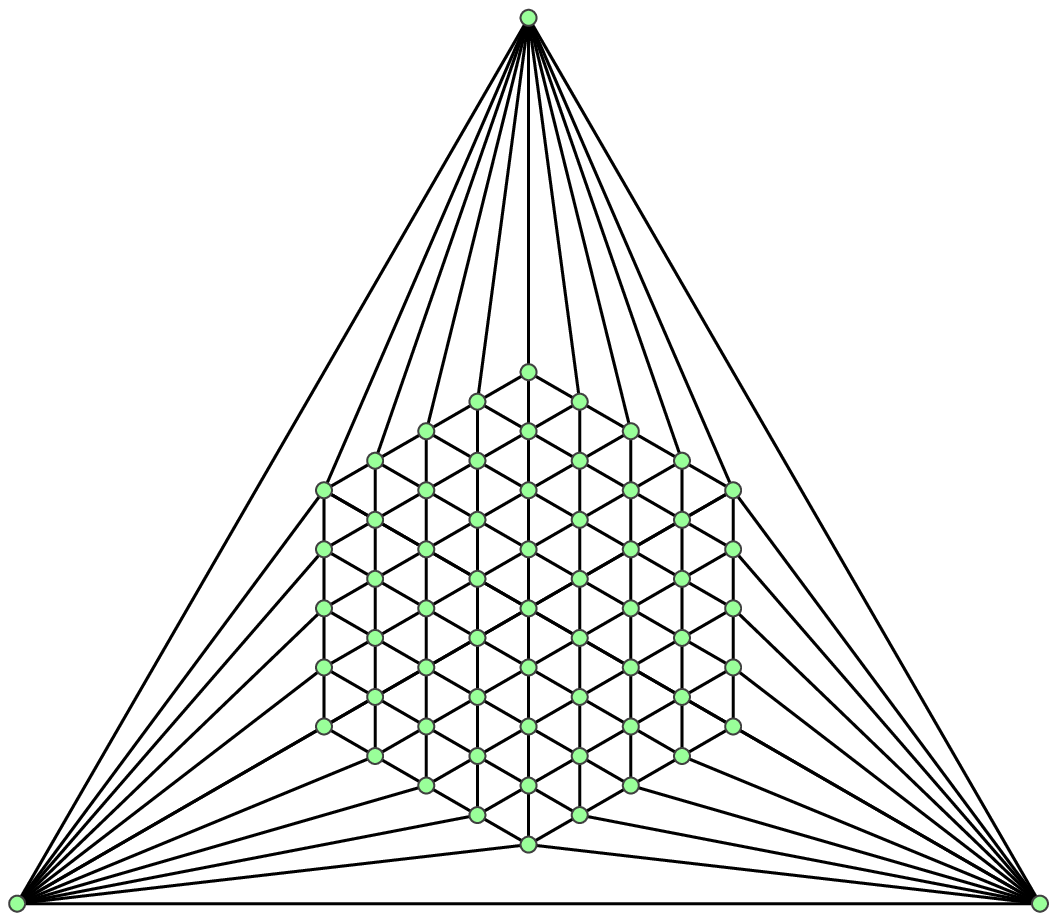}}{A plane triangulation with only \Oh{\sqrt{n}} segments.}


We now prove that the bound on the number of slopes in \thmref{ThreeConnectedSegmentsSlopes} is tight up to an additive constant.

\begin{lemma}
\lemlabel{TriangulationsSlopes}
For all $n\geq3$, there is an $n$-vertex  planar triangulation $G_n$ that has at least $n+2$ slopes in every plane drawing. For a particular choice of outerface, there are at least $2n-2$ slopes in every plane drawing. 
\end{lemma}

\begin{proof} Let $G_n$ be the graph with vertex set $\{v_1,v_2,\dots,v_n\}$ and edge set $\{v_1v_i,v_2v_i:3\leq i\leq n\}\cup\{v_iv_{i+1}:1\leq i\leq n-1\}$. $G_n$ is a planar triangulation.  Every $3$-cycle in $G_n$ contains $v_1$ or $v_2$. Thus $v_1$ or $v_2$ is in the boundary of the outerface in every plane drawing of $G_n$.  By \obsref{*}, the number of slopes in any plane drawing  of $G_n$ is at least $(n-3)+1+1+3=n+2$. As illustrated in \figref{Fish}(a), if we fix the outerface of $G_n$ to be $(v_1,v_2,v_n)$, then  the number of slopes is at least  $(n-3)+(n-3)+1+3=2n-2$ slopes by \obsref{*} 
\end{proof}

As illustrated in \figref{Fish}(b), the graph $G_n$ in \lemref{TriangulationsSlopes} has a plane drawing (using a different embedding) with only $\ceil{\frac{3n}{2}}$ slopes. 

\Figure{Fish}{\includegraphics{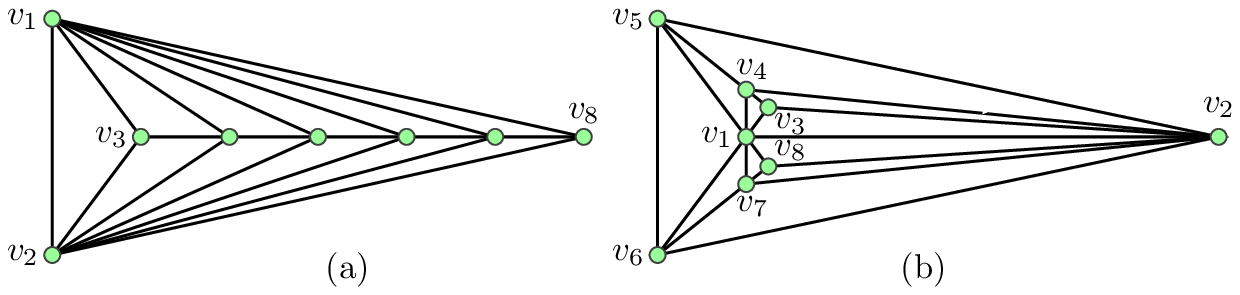}}{The graph $G_8$ from \lemref{TriangulationsSlopes}.}

Since deleting an edge from a drawing cannot increase the number of slopes, and every plane graph can be triangulated to a $3$-connected plane graph, \thmref{ThreeConnectedSegmentsSlopes} implies:

\begin{corollary} 
Every $n$-vertex plane graph has a plane drawing  with at most $2n-10$ slopes.\qed
\end{corollary}

\begin{open}
Is there some $\epsilon>0$, such that  every $n$-vertex planar triangulation
has a plane drawing with $(2-\epsilon)n+\Oh{1}$ slopes?
\end{open}

On the other hand, \thmref{ThreeConnectedSegmentsSlopes} does not imply any upper bound on the  number of segments for all planar graphs.  A natural question to ask is whether \thmref{ThreeConnectedSegmentsSlopes} can be extended to plane graphs that are not $3$-connected. We have the following lower bound.

\begin{lemma}
\lemlabel{BadTwoConnected}
For all even $n\geq 4$, there is a $2$-connected plane graph with $n$ vertices (and $\frac{5}{2}n-4$ edges) that has as many segments as edges in every drawing.
\end{lemma}

\begin{proof} Let $G_n$ be the graph with vertex set $\{v,w,x_i,y_i:1\leq i\leq\half(n-2)\}$ and edge set $\{vw,x_iy_i,vx_i,vy_i,wx_i,wy_i:1\leq i\leq\half(n-2)\}$. Consider the plane embedding of $G_n$ with the cycle $(v,w,y_n)$ as the outerface, as illustrated in \figref{BadTwoConnected}.  Since the outerface is a triangle, no two edges incident to $v$ can share a segment, and  no two edges incident to $w$ can share a segment. Consider two edges $e$ and $f$ both incident to a vertex $x_i$ or $y_i$. The endpoints of $e$ and $f$ induce a triangle. Thus $e$ and $f$ cannot share a segment. Therefore no two edges in $G_n$ share  a segment. \end{proof}

\Figure{BadTwoConnected}{\includegraphics{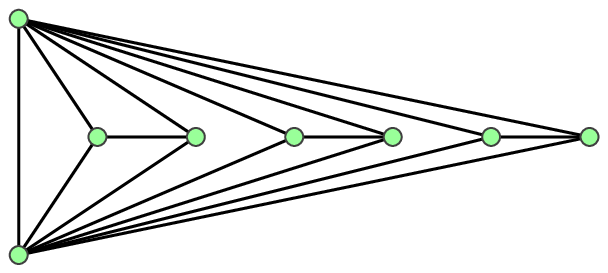}}{The graph $G_8$ in \lemref{BadTwoConnected}.}

Note that the drawing technique from \figref{TwoTree} can be used to draw the
graph $G_n$ in \lemref{BadTwoConnected} with only $2n+\Oh{1}$  segments.

\begin{open} 
What is the minimum $c$ such that every $n$-vertex plane (or planar) graph has
a plane drawing with at most $cn+\Oh{1}$ segments?
\end{open}

\mySubSection{Cubic $3$-Connected Plane Graphs}{Cubic}

A graph in which every vertex has degree three is \emph{cubic}.  It is easily seen that \thmref{ThreeConnectedS} implies that every cubic plane $3$-connected graph on $n$ vertices has a plane drawing with at most $\frac{5}{4}n+\Oh{1}$ segments. This result can be improved as follows.

%

\begin{lemma}
\lemlabel{CubicThreeConnectedS}
Every cubic plane $3$-connected graph $G$ on $n$ vertices has a plane
drawing with at most $n+2$ segments.
\end{lemma}

\begin{proof} Let $D$ be the plane drawing of $G$ from \thmref{ThreeConnectedS}. Recall the definitions and the arguments for counting segments in \thmref{ThreeConnectedS}. By \eqnref{seg}, the number of segments is at most
\begin{equation*}
1+h_s+\sum_{i=2}^K (|\pred{V_i}|-1)\enspace.
\end{equation*}  
By the properties of the canonical ordering for plane cubic graphs, $|\pred{V_i}|=2$ for all $2\leq i\leq K-1$, and $|\pred{V_K}|=3$. Thus $|R|\leq 1$. As in \thmref{ThreeConnectedS}, the number of segment-heavy steps with $|\pred{V_i}|=2$ is at most $|S|+1$. Thus $h_s\leq|S|+2$. Therefore the number of segments in $D$ is at most
\begin{equation*}
1+(|S|+2) + (K-2) + 2\;=\;|S|+3+K \;=\;
|S|+3+n-1-|S| \;=\; n+2\enspace,
\end{equation*}
as claimed.
\end{proof}

Our bound on the number of slopes in a drawing of a $3$-connected plane graph (\thmref{ThreeConnectedSegmentsSlopes}) can be drastically improved when the graph is cubic.

\begin{theorem}
\thmlabel{CubicThreeConnectedPlanarSlopes}
Every cubic $3$-connected plane graph has a plane drawing in which every edge has slope in $\{\frac{\pi}{4},\frac{\pi}{2},\frac{3\pi}{4}\}$, except for three edges on the outerface. 
\end{theorem}

\begin{proof} Let $\sigma=(V_1,V_2,\dots,V_K)$ be a canonical ordering of $G$. We re-use the notation from \thmref{ThreeConnectedS}, except that a representative vertex of $V_i$ may be the first or last vertex in $V_i$. Since $G$ is cubic, $|P(V_i)|=2$ for all $1<i<K$, and every vertex not in $\{v_1,v_2,v_n\}$ has exactly one successor. We proceed by induction on $i$ with the hypothesis that $G_i$ has a plane drawing $D_i$ that satisfies the following invariants.

\begin{description}

\item[Invariant 1:] $C_i\setminus v_1v_2$ is \emph{$X$-monotone}; that is, the path from $v_1$ to $v_2$ in $C_i\setminus v_1v_2$  has non-decreasing X-coordinates.

\item[Invariant 2:] Every peak in $D_i$, $i<K$, has a successor.

\item[Invariant 3:] If there is a vertical edge above $v$ in $D_i$, then  all the edges of $G$ that are incident to $v$ are in $G_i$. 

\item[Invariant 4:] $D_i$ has no edge crossings.

\end{description}

Let $D_2$ be the drawing of $G_2$ constructed as follows. Draw $v_1v_2$ horizontally with $X(v_1)<X(v_2)$. This accounts for one edge whose slope is not in $\{\frac{\pi}{4},\frac{\pi}{2},\frac{3\pi}{4}\}$.  Now draw $v_1v_3$ with slope $\frac{\pi}{4}$, and draw $v_2v_3$ with slope $\frac{3\pi}{4}$. Add any division vertices on the segment $v_1v_3$. Now $v_3$ is the only peak in $D_2$, and it has a successor by the definition of the canonical ordering. Thus all the invariants are satisfied for the base case $D_2$.

Now suppose that $2<i<K$ and we have a drawing of $D_{i-1}$ that satisfies the invariants. Suppose that $\pred{V_i}=\{u,w\}$, where $u$ and $w$ are the left and the right predecessors of $V_i$, respectively. Without loss of generality,  $Y(w)\leq Y(u)$. Let the representative vertex $v_i$ be last vertex in $V_i$. Position $v_i$ at the intersection of a vertical segment above $w$, and a segment of slope $\frac{\pi}{4}$ from $u$, and add any division vertices on $\overline{uv_i}$, as illustrated in \figref{ThreeRegular}(a). Note that there is no vertical edge above $w$ by invariant~3 for $D_{i-1}$. (For the case in which $Y(u)<Y(w)$, we take the representative vertex $v_i$ to be the first vertex in $V_i$, and the edge $wv_i$ has slope $\frac{3\pi}{4}$, as illustrated in \figref{ThreeRegular}(b).)\ 

\Figure{ThreeRegular}{\includegraphics{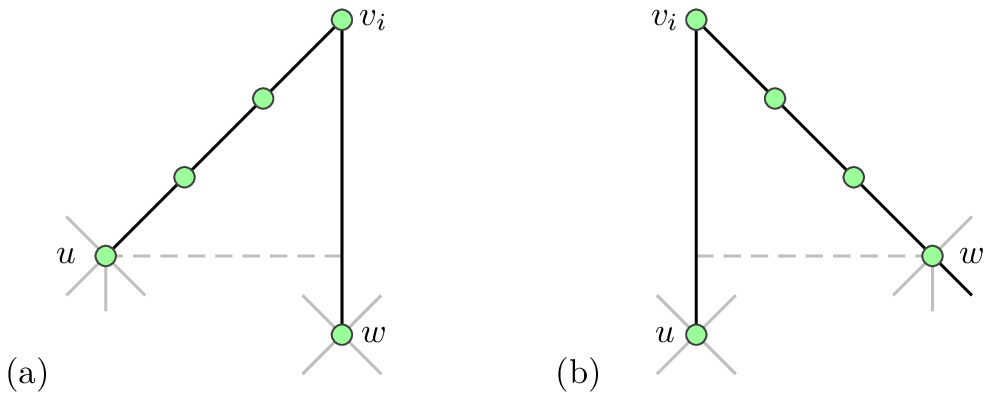}}{Construction of a $3$-slope drawing of a cubic $3$-connected plane graph.}

Clearly the resulting drawing $D_i$ is $X$-monotone. Thus invariant~1 is maintained.  The vertex $v_i$ is the only peak in $D_i$ that is not a peak in $D_{i-1}$. Since $v_i$ has a successor by the definition of the canonical ordering, invariant~2 is maintained. The vertical edge $wv_i$ satisfies invariant~3, since $v_i$ is the sole successor of $w$. Thus invariant~3 is maintained. No vertex between $u$ and $w$ (on the path from $u$ to $w$ in $C_{i-1}\setminus v_1v_2$) is higher than the higher of $u$ and $w$. Otherwise there would be a peak, not equal to $v_n$, with no successor, and thus violating invariant~2 for $D_{i-1}$. Thus the edges in $D_i\setminus D_{i-1}$ do not cross any edges in $D_i$. In particular, there is no edge $ux$ in $D_{i-1}$ with slope $\frac{\pi}{4}$ and $Y(x)>Y(u)$. 

It remains to draw the vertex $v_n$. Suppose $v_n$ is adjacent to $v_1$, $u$, and $w$, where $X(v_1)<X(u)<X(w)$.  By invariants~1 and 3 applied to $v_1$, $u$ and $w$, there is point $p$ vertically above $u$ that is visible from $v_1$ and $w$. Position $v_n$ at $p$ and draw its incident edges. We obtain the desired drawing of $G$. The edge $v_nu$ has slope $\frac{\pi}{2}$, while $v_nv_1$ and $v_nw$ are the remaining two edges whose slope is not in $\{\frac{\pi}{4},\frac{\pi}{2},\frac{3\pi}{4}\}$. \end{proof}


A number of notes regarding \thmref{CubicThreeConnectedPlanarSlopes} are in order:

\begin{itemize}

\item By \lemref{ThreeSlope} we could have used any set of three slopes instead of $\{\frac{\pi}{4},\frac{\pi}{2},\frac{3\pi}{4}\}$ in \thmref{CubicThreeConnectedPlanarSlopes}. 

\item By \obsref{*}, the bound of six on the number of slopes in \thmref{CubicThreeConnectedPlanarSlopes} is optimal for any $3$-connected cubic plane graph whose outerface is a triangle. It is easily seen that there is such a graph on $n$ vertices for all even $n\geq4$.

\item \thmref{CubicThreeConnectedPlanarSlopes} was independently obtained by \citet{Kant-WG92}. We believe that our proof is much simpler. \citet{Kant-WG92} also claimed to prove that every plane graph with maximum degree three has a $3$-slope drawing (except for one bent edge). This claim is false. Consider the plane graph $G$ illustrated in \figref{BadKant}(a). It is easily seen that $G$ has no $3$-slope plane (straight-line) drawing. Thus the cubic plane graph illustrated in \figref{BadKant}(b), which contains a linear number of copies of $G$, must have a linear number of bends in any plane drawing on three slopes.

\citet{Kant-WG92} also claimed to prove that every planar graph with maximum degree three (except $K_4$) has a drawing in which every angle (between consecutive edges incident to a vertex) is at least $\frac{\pi}{3}$, except for at most four angles. The example in \figref{BadKant}(b) is a counterexample to this claim as well. It is easily seen that every drawing of $G$ has an angle less than $\frac{\pi}{3}$. (Assume otherwise, and start from back-to-back drawings of two equilateral triangles.)\ Thus the cubic plane graph illustrated in \figref{BadKant}(b) has a linear number of angles less than $\frac{\pi}{3}$. 

\end{itemize}

\Figure{BadKant}{\includegraphics{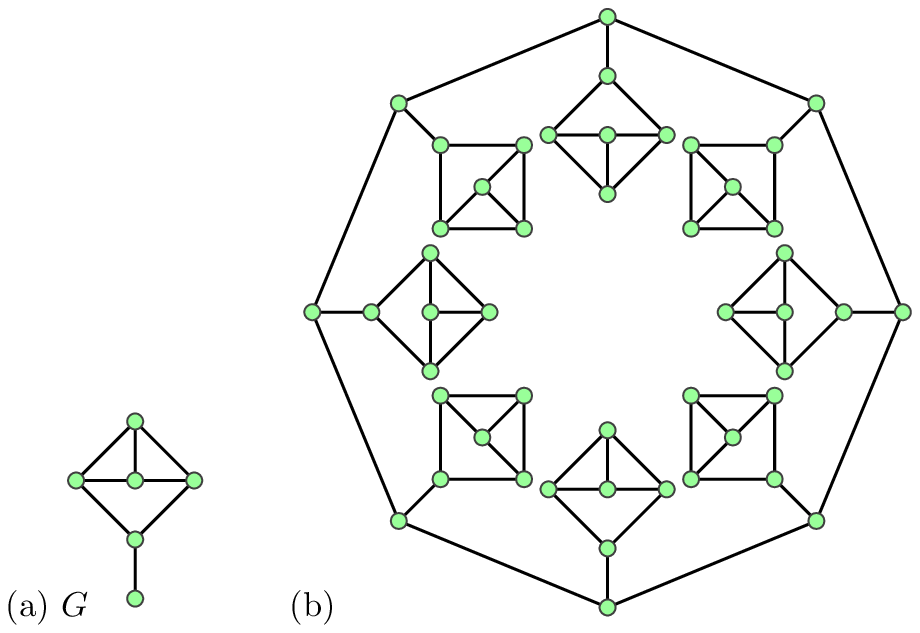}}{Counterexample to the claim of Kant \citep{Kant-WG92}.}


\begin{corollary}
\corlabel{CubicBends}
Every cubic $3$-connected plane graph has a plane `drawing' with three slopes and three bends on the outerface. 
\end{corollary}

\begin{proof} Apply the proof of \thmref{CubicThreeConnectedPlanarSlopes} with two exceptions. First the edge $v_1v_2$ is drawn with one bend. The segment incident to $v_1$ has slope $\frac{3\pi}{4}$, and the segment incident to $v_2$ has slope $\frac{\pi}{4}$. The second exception regards how to draw the edges incident to $v_n$. Suppose $v_n$ is adjacent to $v_1$, $u$, and $w$, where $X(v_1)<X(u)<X(w)$. There is a point $s$ above $v_1$, a point $p$ above $u$, and a point $t$ above $w$, so that the slope of $sp$ is $\frac{\pi}{4}$ and the slope of $tp$ is $\frac{3\pi}{4}$. Place $v_n$ at $p$, draw the edge $v_nu$ vertical, draw the edge $v_1v_n$ with one bend through $s$ (with slopes $\{\frac{\pi}{2},\frac{\pi}{4}\}$), and draw the edge $wv_n$ with one bend through $t$ (with slopes $\{\frac{\pi}{2},\frac{3\pi}{4}\}$). \end{proof}

\begin{open}
Does there exist a function $f$ such that every plane graph with maximum
degree $\Delta$ has a plane drawing with $f(\Delta)$ slopes? This is open even
for maximal outerplanar graphs. Note that there exist bounded degree (non-planar) graphs for which the number of slopes is unbounded in every drawing \citep{PachPal-EJC06,BMW-EJC06,DSW-Slopes}. The best bounds are in our companion paper \citep{DSW-Slopes}, in which we prove that there exists $\Delta$-regular $n$-vertex graphs with at least $n^{1-\frac{8+\epsilon}{\Delta+4}}$ slopes in every drawing.
\end{open}

\begin{open}
In all our results, we have not studied other aesthetic criteria such as symmetry and small area (with the vertices at grid points). Many open problems remain when combining ``few slopes or segments'' with other aesthetic criteria. For example, can \thmref{CubicThreeConnectedPlanarSlopes} be generalised to prove that every cubic $3$-connected plane graph on $n$ vertices has a plane grid drawing with polynomial (in $n$) area, such that every edge has one of three slopes (except for three edges on the outerface)?
\end{open}

\section*{Acknowledgements} 

This research was initiated at the \emph{International Workshop on Fixed Parameter Tractability in Geometry and Games}, organised by Sue Whitesides; Bellairs Research Institute of McGill University, Barbados, February 7--13, 2004. Thanks to all of the participants for creating a stimulating working environment. Special thanks to  Mike Fellows for suggesting the problem. Thanks to Therese Biedl for pointing out reference \citep{Kant-WG92}. 


\def\soft#1{\leavevmode\setbox0=\hbox{h}\dimen7=\ht0\advance \dimen7
  by-1ex\relax\if t#1\relax\rlap{\raise.6\dimen7
  \hbox{\kern.3ex\char'47}}#1\relax\else\if T#1\relax
  \rlap{\raise.5\dimen7\hbox{\kern1.3ex\char'47}}#1\relax \else\if
  d#1\relax\rlap{\raise.5\dimen7\hbox{\kern.9ex \char'47}}#1\relax\else\if
  D#1\relax\rlap{\raise.5\dimen7 \hbox{\kern1.4ex\char'47}}#1\relax\else\if
  l#1\relax \rlap{\raise.5\dimen7\hbox{\kern.4ex\char'47}}#1\relax \else\if
  L#1\relax\rlap{\raise.5\dimen7\hbox{\kern.7ex
  \char'47}}#1\relax\else\message{accent \string\soft \space #1 not
  defined!}#1\relax\fi\fi\fi\fi\fi\fi} \def\cprime{$'$}

\end{document}

%% file: myTheorems.tex
\usepackage{amsthm}

\theoremstyle{plain}

\newtheorem{observation}{Observation}
\newtheorem{theorem}{Theorem}

\newtheorem{claim}{Claim}
\newtheorem{lemma}{Lemma}
\newtheorem{corollary}{Corollary}

\theoremstyle{definition}
\newtheorem{conjecture}{Conjecture}

\newtheorem{open}{Open Problem}

%% file: myLabels.tex
\newcommand{\thmlabel}[1]{\label{thm:#1}}
\newcommand{\thmref}[1]{Theorem~\ref{thm:#1}}

\newcommand{\lemlabel}[1]{\label{lem:#1}}
\newcommand{\lemref}[1]{Lemma~\ref{lem:#1}}

\newcommand{\eqnlabel}[1]{\label{eqn:#1}}
\newcommand{\eqnref}[1]{\eqref{eqn:#1}}

\newcommand{\figlabel}[1]{\label{fig:#1}}
\newcommand{\figref}[1]{Figure~\ref{fig:#1}}

\newcommand{\seclabel}[1]{\label{sec:#1}}
\newcommand{\secref}[1]{Section~\ref{sec:#1}}

\newcommand{\mySection}[2]{\section{\boldmath\textsc{#1}}\seclabel{#2}}
\renewcommand{\mySection}[2]{\section{#1}\seclabel{#2}}
\newcommand{\mySubSection}[2]{\subsection{\boldmath #1}\seclabel{#2}}

\newcommand{\corlabel}[1]{\label{cor:#1}}
\newcommand{\corref}[1]{Corollary~\ref{cor:#1}}

\newcommand{\tablabel}[1]{\label{tab:#1}}
\newcommand{\tabref}[1]{Table~\ref{tab:#1}}

\newcommand{\obslabel}[1]{\label{obs:#1}}
\newcommand{\obsref}[1]{Observation~\ref{obs:#1}}

%% file: myCommands.tex
\DeclareFontShape{OT1}{cmr}{b}{sc}
{<5><6><7><8><9><10><12><10.95><14.4><17.28><20.74><24.88>cmbcsc10}{}
\DeclareFontShape{OT1}{cmr}{bx}{sc}
{<5><6><7><8><9><10><12><10.95><14.4><17.28><20.74><24.88>cmbcsc10}{}

\DeclareFontShape{OT1}{cmr}{b}{tt}
{<5><6><7><8><9><10><12><10.95><14.4><17.28><20.74><24.88>cmbtt10}{}
\DeclareFontShape{OT1}{cmr}{bx}{tt}
{<5><6><7><8><9><10><12><10.95><14.4><17.28><20.74><24.88>cmbtt10}{}

\newcommand{\NP}{\ensuremath{\mathcal{NP}}}

\newcommand{\FancyFootnotes}{\renewcommand{\thefootnote}{\fnsymbol{footnote}}}
\newcommand{\FancyFootnotesOff}{\renewcommand{\thefootnote}{\arabic{footnote}}}

\newcommand{\paran}[1]{\textup{(}#1\textup{)}}

\newcommand{\ceil}[1]{\ensuremath{\protect\lceil#1\rceil}}

\newcommand{\floor}[1]{\ensuremath{\protect\lfloor#1\rfloor}}

\newcommand{\half}{\ensuremath{\protect\tfrac{1}{2}}}

\newcommand{\Oh}[1]{\ensuremath{\protect\mathcal{O}(#1)}}

\newcommand{\Figure}[4][htb]{
\begin{figure}[#1]
	\vspace*{1ex}
	\begin{center}#3\end{center}
	\vspace*{-1ex}
	\caption{\figlabel{#2}#4}
\end{figure}
}

\newlength{\marginboxwidth}

\usepackage{pifont}

\newcommand{\COMMENT}[1]{\medskip\noindent\framebox{\parbox{\textwidth-2mm}{#1}}\medskip}


%% file: DESW-SlopesSegments.bbl
\begin{thebibliography}{26}
\providecommand{\natexlab}[1]{#1}
\providecommand{\url}[1]{\texttt{#1}}
\providecommand{\urlprefix}{}
\expandafter\ifx\csname urlstyle\endcsname\relax
  \providecommand{\doi}[1]{doi:\discretionary{}{}{}#1}\else
  \providecommand{\doi}{doi:\discretionary{}{}{}\begingroup
  \urlstyle{rm}\Url}\fi

\bibitem[{Bar{\'a}t et~al.(2006)Bar{\'a}t, Matou{\v{s}}ek, and
  Wood}]{BMW-EJC06}
\textsc{J{\'a}nos Bar{\'a}t, Ji{\v{r}}{\'i} Matou{\v{s}}ek, and David~R. Wood}.
\newblock Bounded-degree graphs have arbitrarily large geometric thickness.
\newblock \emph{Electron. J. Combin.}, 13(1):R3, 2006.

\bibitem[{Bhasker and Sahni(1987)}]{BS-Networks87}
\textsc{Jayaram Bhasker and Sartaj Sahni}.
\newblock A linear time algorithm to check for the existence of a rectangular
  dual of a planar triangulated graph.
\newblock \emph{Networks}, 17(3):307--317, 1987.

\bibitem[{Bhasker and Sahni(1988)}]{BS-Algo88}
\textsc{Jayaram Bhasker and Sartaj Sahni}.
\newblock A linear algorithm to find a rectangular dual of a planar
  triangulated graph.
\newblock \emph{Algorithmica}, 3(2):247--278, 1988.

\bibitem[{Czyzowicz et~al.(1998)Czyzowicz, Kranakis, and Urrutia}]{CKU-IPL98}
\textsc{Jurek Czyzowicz, Evangelos Kranakis, and Jorge Urrutia}.
\newblock A simple proof of the representation of bipartite planar graphs as
  the contact graphs of orthogonal straight line segments.
\newblock \emph{Inform. Process. Lett.}, 66(3):125--126, 1998.

\bibitem[{de~Castro et~al.(2002)de~Castro, Cobos, Dana, M{\'a}rquez, and
  Noy}]{CCDMN-JGAA02}
\textsc{Natalia de~Castro, Francisco~Javier Cobos, Juan~Carlos Dana, Alberto
  M{\'a}rquez, and Marc Noy}.
\newblock Triangle-free planar graphs as segment intersection graphs.
\newblock \emph{J. Graph Algorithms Appl.}, 6(1):7--26, 2002.

\bibitem[{de~Fraysseix and {Ossona de Mendez}(2004)}]{dFdM-GD04}
\textsc{Hubert de~Fraysseix and Patrice {Ossona de Mendez}}.
\newblock Contact and intersection representations.
\newblock In \textsc{J\'{a}nos Pach}, ed., \emph{Proc. 12th International Symp.
  on Graph Drawing (GD '04)}, vol. 3383 of \emph{Lecture Notes in Comput.
  Sci.}, pp. 217--227. Springer, 2004.

\bibitem[{de~Fraysseix et~al.(1995)de~Fraysseix, {Ossona de Mendez}, and
  Pach}]{dFdMP95}
\textsc{Hubert de~Fraysseix, Patrice {Ossona de Mendez}, and J\'{a}nos Pach}.
\newblock A left-first search algorithm for planar graphs.
\newblock \emph{Discrete Comput. Geom.}, 13(3-4):459--468, 1995.

\bibitem[{{de Fraysseix} et~al.(1990){de Fraysseix}, Pach, and
  Pollack}]{dFPP90}
\textsc{Hubert {de Fraysseix}, J\'{a}nos Pach, and Richard Pollack}.
\newblock How to draw a planar graph on a grid.
\newblock \emph{Combinatorica}, 10(1):41--51, 1990.

\bibitem[{Dirac(1961)}]{Dirac61}
\textsc{Gabriel~A. Dirac}.
\newblock On rigid circuit graphs.
\newblock \emph{Abh. Math. Sem. Univ. Hamburg}, 25:71--76, 1961.

\bibitem[{Dujmovi{\'c} et~al.(2005)Dujmovi{\'c}, Suderman, and
  Wood}]{DSW-Slopes}
\textsc{Vida Dujmovi{\'c}, Matthew Suderman, and David~R. Wood}.
\newblock Graph drawings with few slopes.
\newblock Submitted, 2005.
\newblock \urlprefix\url{http://arxiv.org/math/0606446}.

\bibitem[{Dujmovi{\'c} and Wood(2006)}]{DujWoo-GD05}
\textsc{Vida Dujmovi{\'c} and David~R. Wood}.
\newblock Graph treewidth and geometric thickness parameters.
\newblock In \textsc{Patrick Healy and Nikola~S. Nikolov}, eds., \emph{Proc.
  13th International Symp. on Graph Drawing (GD '05)}, vol. 3843 of
  \emph{Lecture Notes in Comput. Sci.}, pp. 129--140. Springer, 2006.
\newblock \urlprefix\url{http://arxiv.org/math/0503553}.

\bibitem[{F{\'a}ry(1948)}]{Fary48}
\textsc{Istv{\'a}n F{\'a}ry}.
\newblock On straight line representation of planar graphs.
\newblock \emph{Acta Univ. Szeged. Sect. Sci. Math.}, 11:229--233, 1948.

\bibitem[{Garg and Tamassia(2001)}]{GT-SJC01}
\textsc{Ashim Garg and Roberto Tamassia}.
\newblock On the computational complexity of upward and rectilinear planarity
  testing.
\newblock \emph{SIAM J. Comput.}, 31(2):601--625, 2001.

\bibitem[{Hartman et~al.(1991)Hartman, Newman, and Ziv}]{HNZ-DM91}
\textsc{Irith Ben-Arroyo Hartman, Ilan Newman, and Ran Ziv}.
\newblock On grid intersection graphs.
\newblock \emph{Discrete Math.}, 87(1):41--52, 1991.

\bibitem[{Kant(1993)}]{Kant-WG92}
\textsc{Goos Kant}.
\newblock Hexagonal grid drawings.
\newblock In \textsc{Ernst~W. Mayr}, ed., \emph{Proc. 18th International
  Workshop in Graph-Theoretic Concepts in Computer Science (WG '92)}, vol. 657
  of \emph{Lecture Notes in Comput. Sci.}, pp. 263--276. Springer, 1993.
\newblock Also see Tech.\ Rep.\ RUU-CS-92-06, Universiteit Utrecht,
  Netherlands, 1992.

\bibitem[{Kant(1996)}]{Kant96}
\textsc{Goos Kant}.
\newblock Drawing planar graphs using the canonical ordering.
\newblock \emph{Algorithmica}, 16(1):4--32, 1996.

\bibitem[{Ossona~de Mendez and de~Fraysseix(1999)}]{OdMdF99}
\textsc{Patrice Ossona~de Mendez and Hubert de~Fraysseix}.
\newblock Intersection graphs of {J}ordan arcs.
\newblock In \emph{Contemporary trends in discrete mathematics}, vol.~49 of
  \emph{DIMACS Ser. Discrete Math. Theoret. Comput. Sci.}, pp. 11--28. Amer.
  Math. Soc., 1999.

\bibitem[{Pach and P{\'a}lv{\"o}lgyi(2006)}]{PachPal-EJC06}
\textsc{J{\'a}nos Pach and D{\"o}m{\"o}t{\"o}r P{\'a}lv{\"o}lgyi}.
\newblock Bounded-degree graphs can have arbitrarily large slope numbers.
\newblock \emph{Electron. J. Combin.}, 13(1):N1, 2006.

\bibitem[{Rahman et~al.(1998)Rahman, Nakano, and Nishizeki}]{RNN-CG98}
\textsc{Md.~Saidur Rahman, {Shin-ichi} Nakano, and Takao Nishizeki}.
\newblock Rectangular grid drawings of plane graphs.
\newblock \emph{Comput. Geom.}, 10(3):203--220, 1998.

\bibitem[{Rahman et~al.(1999)Rahman, Nakano, and Nishizeki}]{RNN-JGAA99}
\textsc{Md.~Saidur Rahman, {Shin-ichi} Nakano, and Takao Nishizeki}.
\newblock A linear algorithm for bend-optimal orthogonal drawings of
  triconnected cubic plane graphs.
\newblock \emph{J. Graph Algorithms Appl.}, 3:31--62, 1999.

\bibitem[{Rahman et~al.(2000)Rahman, Nakano, and Nishizeki}]{RNN-JAlg00}
\textsc{Md.~Saidur Rahman, Shin-ichi Nakano, and Takao Nishizeki}.
\newblock Box-rectangular drawings of plane graphs.
\newblock \emph{J. Algorithms}, 37(2):363--398, 2000.

\bibitem[{Rahman et~al.(2002)Rahman, Nakano, and Nishizeki}]{RNN-CG02}
\textsc{Md.~Saidur Rahman, {Shin-ichi} Nakano, and Takao Nishizeki}.
\newblock Rectangular drawings of plane graphs without designated corners.
\newblock \emph{Comput. Geom.}, 21(3):121--138, 2002.

\bibitem[{Rahman et~al.(2004)Rahman, Nishizeki, and Ghosh}]{RNG-JAlg04}
\textsc{Md.~Saidur Rahman, Takao Nishizeki, and Shubhashis Ghosh}.
\newblock Rectangular drawings of planar graphs.
\newblock \emph{J. Algorithms}, 50:62--78, 2004.

\bibitem[{Thomassen(1984)}]{Thomassen84}
\textsc{Carsten Thomassen}.
\newblock Plane representations of graphs.
\newblock In \emph{Progress in graph theory}, pp. 43--69. Academic Press,
  Toronto, 1984.

\bibitem[{Ungar(1953)}]{Ungar-JLMS53}
\textsc{Peter Ungar}.
\newblock On diagrams representing maps.
\newblock \emph{J. London Math. Soc.}, 28:336--342, 1953.

\bibitem[{Wagner(1936)}]{Wagner36}
\textsc{Klaus Wagner}.
\newblock Bemerkung zum {V}ierfarbenproblem.
\newblock \emph{Jber. Deutsch. Math.-Verein.}, 46:26--32, 1936.

\end{thebibliography}
